\documentclass[11 pt,a4paper,leqno]{article}
\usepackage{latexsym,amsmath,amscd,amsfonts}

\RequirePackage{amsmath}
\RequirePackage{amssymb}
\RequirePackage{ifthen}
\RequirePackage{theorem}
\RequirePackage{pslatex}

\usepackage[all]{xy}

\title{\'Equations aux $q$-diff\'erences lin\'eaires:
factorisation, r\'esolution et th\'eor\`emes d'indices}

\author{\small Jacques Sauloy, \\
\small Laboratoire Emile Picard, UMR 5580,\\
\small Universit\'e Paul Sabatier, U.F.R. M.I.G.\\
\small 118, route de Narbonne, 31062 Toulouse Cedex 4, France\\
\small {\tt sauloy@math.univ-toulouse.fr}
}

\date{10 d\'ecembre 2008}

\theorembodyfont{\sl} 
\newtheorem{thm}{Th\'eor\`eme}[section]
\newtheorem{lemma}[thm]{Lemme}
\newtheorem{prop}[thm]{Proposition}
\newtheorem{cor}[thm]{Corollaire}

\theorembodyfont{\rm}
\newtheorem{rmq}[thm]{Remarque}

\newtheorem{defn}[thm]{D\'efinition}

\numberwithin{equation}{thm}

\def \Pr {\textsl{Preuve. - }}

\def\sq{\sigma_q}

\def\C{{\mathbf C}}
\def\Q{{\mathbf Q}}
\def\Z{{\mathbf Z}}
\def\R{{\mathbf R}}
\def\N{{\mathbf N}}

\def\E{{\mathbf E}}
\def\Eq{{\E_{q}}}
\def\M{{\mathcal M}}
\def\O{{\mathcal O}}
\def\D{{\mathcal D}}

\def\tq{\, | \,}
\def\lmod{\left |}
\def\rmod{\right |}

\def\Coker{{\text{Coker~}}}

\newcommand{\Ker}{\text{Ker }}

\begin{document}

\maketitle

\bigskip \hrule \bigskip

\centerline{\textbf{\emph{R\'esum\'e}}}

\emph{Nous d\'ecrivons des algorithmes explicites pour la factorisation
d'op\'erateurs et la r\'esolution d'\'equations aux $q$-diff\'erences.
Il s'agit d'une pr\'esentation ``concr\`ete'' de r\'esultats bien connus.}

\bigskip \hrule \bigskip

\centerline{\textbf{\emph{Abstract}}}

\emph{We describe explicit algorithms for factoring $q$-difference
operators and solving $q$-difference equations. These are well known
results, presented in a ``concrete'' form.}

\bigskip \hrule \bigskip

\tableofcontents

\bigskip \hrule \bigskip

%%%%%%%%%%%%%%%%%%%%%%%%%%%%%%%%%%%%%%%%%%%%%%%%%%%%%%%%%%%%%%%%%%%%%%%%%%%%%%

% 1

\section{Introduction et g\'en\'eralit\'es}

% 1.1

\subsection{Introduction}

Une \'equation aux $q$-diff\'erences lin\'eaire d'ordre $n$
est une \'equation fonctionnelle de la forme:

\begin{equation} 
\label{equation:eqd}
a_{n}(z) f(q^{n} z) + a_{n-1}(z) f(q^{n-1} z) + \cdots + a_{0}(z) f(z) = 0.
\end{equation}

Il y a plusieurs mani\`eres de justifier leur \'etude. La plus ancienne
met en jeu le ``$q$-calcul'', tout particuli\`erement des formules
d\'ecouvertes par Euler dans l'\'etude des partitions, donc \`a la
fronti\`ere de la combinatoire et de l'arithm\'etique. Si ces formules
elles-m\^emes ne font en g\'en\'eral apparaitre qu'un param\`etre,
not\'e $q$, certaines s\'eries g\'en\'eratrices qui leurs sont associ\'ees
satisfont des \'equations fonctionnelles du type (\ref{equation:eqd})
l\`a o\`u les s\'eries g\'en\'eratrices construites \`a partir de
coefficients combinatoires plus classiques satisferaient des \'equations
diff\'erentielles. \\

A titre d'exemple, un ``$q$-analogue'' de la factorielle
$n ! = \underset{1 \leq i \leq n}{\prod} i$ est la ``$q$-factorielle''
$[n]_{q} !: = \underset{1 \leq i \leq n}{\prod} \frac{q^{i}-1}{q-1}$.
La s\'erie $\underset{n \geq 0}{\sum} \frac{z^{n}}{n!}$ a pour somme 
la fonction $f(z) = e^{z}$ qui est solution de l'\'equation diff\'erentielle 
$f' = f$.
De m\^eme, la s\'erie $\underset{n \geq 0}{\sum} \frac{z^{n}}{[n]_{q}!}$
admet pour somme une fonction $f$, l'un des (nombreux) $q$-analogues
de l'exponentielle, et celle-ci est solution de l'\'equation fonctionnelle
$\frac{f(qz) - f(z)}{(q-1)z} = f(z)$, qui se ram\`ene facilement \`a
la forme (\ref{equation:eqd}). Si l'on fait tendre $q$ vers $1$ dans la
s\'erie ou dans l'\'equation fonctionnelle, on retrouve sa contrepartie
classique. De nombreuses autres motivations existent, plus s\'erieuses 
ou aussi amusantes, plus modernes ou aussi classiques, plus conceptuelles 
ou aussi ph\'enom\'enologiques. On en trouvera un tableau introductif,
avec un survol des grands probl\`emes dans ce domaine, dans l'article
\cite{DRSZ}. \\

Dans le pr\'esent article, nous pr\'esenterons de mani\`ere
d\'etaill\'ee quelques algorithmes de r\'esolution explicite d'\'equations
aux $q$-diff\'erences lin\'eaires. On est souvent guid\'e par l'analogie 
avec les \'equations diff\'erentielles, et, comme pour celles-ci,
un outil important est la factorisation d'op\'erateurs, que nous \'etudions
donc aussi. Notre point de vue est analytique et non seulement formel. 
Nous \'etudions ces \'equations dans le champ complexe et cherchons 
des solutions holomorphes (dans certains domaines). On a constat\'e 
depuis longtemps que presque rien n'\'etait 
possible\footnote{Une perc\'ee a cependant \'et\'e r\'ecemment 
effectu\'ee par Lucia Di Vizio dans \cite{LDVnoncom}.}
sans l'hypoth\`ese que $q \in \C^{*}$ est de module diff\'erent de $1$. \\

\textbf{Exercice. -} D\'eterminer toutes les fonctions $f$ m\'eromorphes
sur un voisinage de $0 \in \C$ telles que $f(qz) = f(z)$ lorsque 
$\lmod q \rmod = 1$. On sera amen\'e \`a traiter \`a part le cas o\`u 
$q$ est une racine de l'unit\'e. \\

Enfin, conform\'ement \`a l'\'evolution de la th\'eorie ``moderne''
des \'equations fonctionnelles lin\'eaires, nous abordons le calcul
des indices: il ne s'agit pas seulement de savoir si les solutions
existent et sont uniques, mais de quantifier le d\'efaut d'existence
ou d'unicité; donc de calculer des dimensions de noyaux et de conoyaux. \\

Bien que nous soyons en permanence guid\'es par la th\'eorie ``classique''
des \'equations diff\'erentielles dans le champ complexe, il y a ici
des sp\'ecificit\'es importantes en comparaison avec le cas classique:

\begin{enumerate}

\item{En supposant seulement les coefficients de (\ref{equation:eqd})
\emph{analytiques} au voisinage de $0$, comme pour les \'equations 
diff\'erentielles, on d\'efinit un polygone de Newton en vue de la
r\'esolution au voisinage de $0$. On rend ses pentes enti\`eres par
un changement de variables $z = w^{\ell}, \ell \in \Z$, puis l'on se 
ram\`ene \`a la recherche de solutions s\'eries par un changement
de fonction inconnue $f = u g$, o\`u la fonction auxiliaire $u$
est solution d'une \'equation \'el\'ementaire. Contrairement au
cas des \'equations diff\'erentielles, nous avons ici la garantie
que \emph{certaines des s\'eries ainsi obtenues seront convergentes}. 
C'est le \emph{lemme d'Adams}, qui a \'et\'e exhum\'e par Marotte et 
Zhang apr\`es un long sommeil et qui a \'et\'e l'une des sources de 
la renaissance de la th\'eorie  analytique (cf. par exemple \cite{RSZ}).
Une version am\'elior\'ee, d\^ue \`a Birkhoff et Guenther, dit que 
tout op\'erateur aux $q$-diff\'erences \`a coefficients analytiques 
admet une factorisation analytique \`a facteurs de degr\'e $1$
(modulo ramification), les coefficients \'etant m\^eme polynomiaux.
Pour les op\'erateurs diff\'erentiels, seule une factorisation formelle 
est en g\'en\'eral possible.}
\item{Supposons les coefficients de (\ref{equation:eqd}) rationnels.
Alors toute solution m\'eromorphe dans un voisinage (\'epoint\'e)
de $0$ dans $\C^{*}$ se prolonge en une unique solution m\'eromorphe
sur tout $\C^{*}$. Ceci tient \`a ce que l'\'equation fonctionnelle
elle-m\^eme permet de prolonger d'un disque $D$ au disque $q D$
ou $q^{-1} D$ (selon que $\lmod q \rmod$ est $> 1$ ou $< 1$). Ainsi, 
l'\'etude locale d'une \'equation aux $q$-diff\'erences comporte 
d'embl\'ee des aspects \emph{globaux} qui n'apparaissent pas dans 
l'\'etude des \'equations diff\'erentielles.}

\end{enumerate}

La plus grande partie des r\'esultats ci-dessous sont classiques:
ils remontent \`a Adams, Birkhoff, Carmichael et d'autres et nous
ne pr\'esentons qu'une mise au propre au go\^ut du jour (XXI\`eme si\`ecle)
d'ailleurs initi\'ee par Bezivin, Ramis, Marotte et Zhang. Il y a
cependant une nouveaut\'e importante de la th\'eorie pr\'esent\'ee
ici par rapport aux travaux plus anciens, c'est que l'on n'y utilise
que des fonctions \emph{uniformes}, ce qui est \'evidemment impossible
pour les \'equations diff\'erentielles. En effet, alors que m\^eme
des \'equations aussi simples que $z f' = 1$ (i.e. $z f'' + f' = 0$)
ou $z f' = \alpha f \;,\; \alpha \not\in \Z$ n'admettent aucune
solution uniforme sur $\C^{*}$ (il faut soit pratiquer des coupures
pour se ramener \`a des ouverts simplement connexes soit passer au
rev\^etement universel), leurs $q$-analogues peuvent \^etre r\'esolues
\`a l'aide de la fonction Theta de Jacobi, selon une suggestion 
de Ramis dans \cite{RamisJPRTraum}. S'il y a une ``monodromie'' 
sous-jacente \`a de telles \'equations, celle-ci devra donc \^etre 
observ\'ee autrement ! C'est l'objet, en particulier, de \cite{JSAIF}, 
\cite{JSGAL} et de la s\'erie \cite{qWGI},\cite{qWGII},\cite{qWGIII}.

\subsection*{Organisation de cet article}

Les sections \ref{subsection:conventions}, \ref{subsection:fncspe}
et \ref{subsection:polygone} de cette introduction contiennent un 
r\'esum\'e rapide des notations, conventions et faits de base; on 
y introduit en particulier les fonctions sp\'eciales uniformes sur 
$\C^{*}$ qui servent, avec les s\'eries enti\`eres, \`a construire toutes
les solutions. Les preuves d\'etaill\'ees figurent dans \cite{JSAIF}
et, en ce qui concerne le polygone de Newton, dans \cite{JSFIL}. \\

Les chapitres \ref{section:factorisation} et \ref{section:resolution} 
sont respectivement consacr\'es \`a la factorisation des op\'erateurs 
et \`a la r\'esolution des \'equations. On \'etudie dans chaque cas
le probl\`eme formel, dont la r\'eponse ressemble \`a celle que l'on
connait pour les \'equations diff\'erentielles, puis le probl\`eme
analytique, et l\`a, les \'enonc\'es sont bien diff\'erents. Pour ces 
deux chapitres, les preuves sont compl\`etes. \\

Enfin, dans le chapitre \ref{section:indices}, on aborde d'un
autre point de vue les \'equations avec second membre: que peut-on
dire du d\'efaut d'existence ou d'unicit\'e de leurs solutions,
autrement dit, du conoyau et du noyau de l'op\'erateur associ\'e ? \\

Aucun de ces r\'esultats n'est tr\`es neuf, j'ai simplement voulu
donner une forme simple, compl\`ete et explicite \`a des \'enonc\'es 
qui apparaissent en filigrane dans \cite{MarotteZhang},dans \cite{JSFIL}
et dans \cite{RSZ}. Aucune connaissance pr\'ealable des \'equations
aux $q$-diff\'erences n'est requise, mais un peu de familiarit\'e avec 
les op\'erateurs diff\'erentiels peut aider \`a suivre les calculs.

\subsection*{Remerciements}

Il y a bien des ann\'ees, Jos\'e-Manuel Aroca et Felipe Cano m'avaient 
invit\'e \`a un expos\'e de propagande pour les $q$-diff\'erences \`a 
Valladolid. Cet article me sert \`a les en remercier. La r\'edaction
finale a \'et\'e achev\'ee lors d'un s\'ejour d'un mois \`a l'Univesit\'e 
des Sciences et Techniques de Lille, o\`u m'avaient invit\'e Changgui
Zhang et Anne Duval, que je remercie pour cela et pour tout le reste.

% 1.2

\subsection{Conventions g\'en\'erales} 
\label{subsection:conventions}

On fixe un nombre complexe $q$ de module $\lmod q \rmod > 1$.
Le corps de base $K$ est l'un des suivants:
$$
\M(\C) \subset \C(\{z\}) \subset \C((z)),
$$
respectivement : fonctions m\'eromorphes sur $\C$, germes de fonctions 
holomorphes en $0$ et s\'eries de Laurent m\'eromorphes formelles sur 
$\C$. Les deux premiers cas constituent le cas ``convergent'',
le dernier cas constitue le cas ``formel''.
Le corps $K$ est muni de la valuation discr\`ete $v_{0}$ (valuation
$z$-adique): on a $v_{0}(f) = k$ si $f = z^{k} g$ avec $g(0) \in K^{*}$. 
On notera $\O$ l'anneau de valuation correspondant, de corps r\'esiduel 
$\O/z \O = \C$. Le corps $K$ est \'egalement muni d'un automorphisme:
$$
\sq : f(z) \mapsto f(qz).
$$
De mani\`ere g\'en\'erale, une extension $(K',\sigma')$ du ``corps
aux diff\'erences'' $(K,\sq)$ est form\'ee d'une extension $K'$ de $K$ 
et d'un automorphisme $\sigma'$ de $K'$ qui \'etend $\sq$. \\

Un premier exemple de telle extension est obtenu par ramification 
de niveau $\ell \in \N^{*}$: on prend $K' := K_{\ell}= K[z_{\ell}]$, 
o\`u $z_{\ell}^{\ell} = z$; et $\sigma'(z_{\ell}) = q_{\ell} z_{\ell}$, 
o\`u $q_{\ell}$ est une racine $\ell$-\`eme de $q$ dans $\C$. 
On peut d'ailleurs rendre ces extensions compatibles, en prenant 
tous les $K_{\ell}$ sous-extensions d'une m\^eme extension $K_{\infty}$
(dans les cas de $\C(\{z\})$ et $\C((z))$, c'est leur cl\^oture
alg\'ebrique), de sorte que l'on ait $z_{ \ell m}^{\ell} = z_{m}$; et
poser $q_{\ell} = e^{r/\ell}$, o\`u l'on a choisi $\tau$ tel que 
$q = e^{r}$, de sorte que $q_{\ell m}^{\ell} = q_{m}$. \\

Un deuxi\`eme exemple de telle extension se pr\'esente lorsque 
l'on choisit une \emph{alg\`ebre} de fonctions $L$ dans laquelle 
on esp\`ere r\'esoudre toutes les \'equations aux $q$-diff\'erences 
\`a coefficients dans $K$. Si $K = \M(\C)$, on prendra ainsi 
$L = \M(\C^{*})$; de m\^eme, si $K= \C(\{z\})$), on prendra 
$L = \M(\C^{*},0)$, le corps des germes en $0$ de fonctions 
m\'eromorphes sur $\C^{*}$. L'automorphisme $\sq$ de $L$ est encore
$f(z) \mapsto f(qz)$. On va construire plus loin des fonctions
dans $L$ solutions des \'equations de base: la fonction Theta de
Jacobi, qui v\'erifie $\sq \Theta = z \Theta$; pour chaque $c \in \C^{*}$, 
une fonction $e_{q,c}$ telle que $\sq e_{q,c} = c e_{q,c}$ (ces fonctions
sont appel\'ees ``$q$-caract\`eres'', en analogie avec les ``caract\`eres''
$z^{\gamma}$); et une fonction $l_{q}$ telle que $\sq l_{q} = l_{q} + 1$
(appel\'ee ``$q$-logarithme''). \\

Si $K = \C((z))$ (``cas formel''), l'extension $L$ de $K$ sera obtenue 
par adjonction de \emph{symboles} permettant de r\'esoudre les \'equations 
de base:
\begin{enumerate}
\item{Pour chaque $c \in \C^{*}$, on introduit $e_{q,c}$ tel que
$\sq e_{q,c} = c e_{q,c}$. On impose de plus les relations suivantes:
$e_{q,c} e_{q,d} = e_{q,cd}$; et $e_{q,1} = 1$, $e_{q,q} = z$.}
\item{On introduit $\Theta$ tel que $\sq \Theta = z \Theta$, et
on le suppose inversible.}
\item{On introduit $l_{q}$ tel que $\sq l_{q} = l_{q} + 1$.}
\end{enumerate}
Dans ce cas, la $K$-alg\`ebre $L$ n'est pas un corps, et m\^eme
pas int\`egre (elle est \'etudi\'ee dans \cite{SVdP}). \\

Les \emph{constantes} de notre th\'eorie sont les \'el\'ements 
invariants par $\sq$. Le corps des constantes $C_{K}$ de $K$ est 
$\C$ dans tous les cas. Celui de $L$ est encore $C_{L} = \C$ dans 
le cas formel (\emph{voir} \cite{SVdP}). Dans le cas convergent, 
il s'identifie au corps $\M(\Eq)$ des fonctions elliptiques
relatif \`a la courbe elliptique $\Eq = \C^{*}/q^{\Z}$
(\emph{voir} \cite{JSAIF}). \\

Tout complexe non nul s'\'ecrit de mani\`ere unique:
$$
c = q^{\epsilon(c)} \overline{c} \quad \text{ avec} \quad
\epsilon(c) \in \Z \quad \text{ et } 
\quad 1 \leq \lmod \overline{c} \rmod < \lmod q \rmod.
$$
Nous identifierons l'ensemble quotient $\Eq = \C^{*}/q^{\Z}$
(qui est une courbe elliptique) \`a la \emph{couronne fondamentale}
$\{z \in \C^{*} \tq 1 \leq \lmod z \rmod < \lmod q \rmod\}$,
qui en est un syst\`eme de repr\'esentants dans $\C^{*}$,
et le repr\'esentant $\overline{c}$ \`a la classe de $c$
modulo $q^{\Z}$. \\

On notera enfin $\D_{q,K} = K\left<\sigma,\sigma^{-1}\right>$ 
l'alg\`ebre de \"{O}re des polyn\^omes de Laurent non commutatifs, 
caract\'eris\'ee par les relations:
$$
\forall x \in K, \forall k \in \Z \;,\;
\sigma^{k} x = \sq^{k}(x) \sigma^{k}.
$$
Un tel polyn\^ome $P \in \D_{q,K}$ mod\'elise donc l'op\'erateur
aux $q$-diff\'erences $P(\sq)$, d'o\`u une op\'eration de $\D_{q,K}$ 
sur $L$. On dira que le polyn\^ome $P$ est \emph{entier} si tous 
ses mon\^omes $\sigma^{k}$ sont \`a degr\'es $k \geq 0$. Pour 
un tel $P = a_{n} \; \sigma^{n} + \cdots + a_{0}$, on a donc:
$$
P.f = P(\sq)(f) = a_{n} \; \sq^{n} f + \cdots + a_{0} \; f = 0,
$$ 
et l'\'equation (\ref{equation:eqd}) s'\'ecrit: $P.f =0$. \\

L'anneau $\D_{q,K}$ est euclidien (\`a gauche et \`a droite). Soit
$P = \underset{\alpha \leq i \leq \beta}{\sum} a_{i} \sigma^{i}$
un \'el\'ement de $\D_{q,K}$. On appellera \emph{degr\'e absolu} 
de $P$ l'entier naturel $\deg(P) = \beta - \alpha$ si 
$a_{\alpha} a_{\beta} \not= 0$ (et $- \infty$ si $P = 0$)
et \emph{valuation $z$-adique} de $P$ l'entier 
$v_{0}(P) = \min(v_{0}(a_{\alpha}),\ldots,v_{0}(a_{\beta}))$
(donc $+ \infty$ si $P = 0$). On a :
$\deg(PQ) = \deg(P) + \deg(Q)$ et $v_{0}(PQ) = v_{0}(P) + v_{0}(Q)$.
(La preuve de toutes ces affirmations est un exercice amusant
laiss\'e au lecteur.)

% 1.3

\subsection{Fonctions sp\'eciales}
\label{subsection:fncspe}

D\'esignons par $\theta_{q}$ la fonction Theta de Jacobi
(\emph{voir} \cite{RamisGrowth}): 
$$
\theta_{q}(z) := \sum_{n \in \mathbf{Z}} (-1)^{n} q^{-n(n-1)/2} z^{n}
$$
C'est une fonction holomorphe sur $\C^{*}$, qui v\'erifie
l'\'equation fonctionnelle:
$$
\theta_{q}(qz) = -qz \theta_{q}(z).
$$
La fonction Theta satisfait la c\'el\`ebre \emph{formule 
du triple produit de Jacobi}:
$$
\theta_{q}(z) = (p;p)_{\infty} (z;p)_{\infty} (pz^{-1};p)_{\infty}
$$
On note ici $p := q^{-1}$ et l'on utilise le \emph{symbole de
Pochhammer} $(x;p)_{\infty} := \prod\limits_{n \geq 0} (1 - p^{n} x)$. 
La fonction $\theta_{q}$ admet donc pour z\'eros les points de
la $q$-spirale logarithmique discr\`ete $q^{\mathbf{Z}}$, et ces
z\'eros sont simples. \\

Modifiant l\'eg\`erement les notations de \cite{JSAIF}, nous poserons
(dans un but de simplicit\'e): 
$$
\Theta_{q}(z) := \theta_{q}(-z/q) = \sum_{n \in \mathbf{Z}} q^{-n(n+1)/2} z^{n}.
$$
C'est une fonction holomorphe sur $\C^{*}$, qui v\'erifie
l'\'equation fonctionnelle:
$$
\Theta_{q}(qz) = z \Theta_{q}(z).
$$
Elle admet pour z\'eros les points de la $q$-spirale logarithmique 
discr\`ete $- q^{\mathbf{Z}}$, et ces z\'eros sont simples. \\

Nous poserons \'egalement: 
$$
l_{q}(z) = z \Theta_{q}'(z)/\Theta_{q}(z),
$$
qui est m\'eromorphe sur $\C^{*}$ avec pour p\^oles (simples) les 
points de la $q$-spirale logarithmique discr\`ete $- q^{\mathbf{Z}}$,
et qui v\'erifie l'\'equation fonctionnelle:
$$
\sq l_{q} = l_{q} + 1.
$$
Enfin, pour tout complexe non nul $c$, nous poserons:
$$
e_{q,c}(z) = \Theta_{q}(z)/\Theta_{q}(c^{-1}z),
$$
qui est m\'eromorphe sur $\C^{*}$ avec pour z\'eros (simples) les 
points de la $q$-spirale logarithmique discr\`ete $- q^{\mathbf{Z}}$,
pour p\^oles (simples) les points de la $q$-spirale logarithmique 
discr\`ete $- c q^{\mathbf{Z}}$, et qui v\'erifie l'\'equation 
fonctionnelle:
$$
\sq e_{q,c} = c e_{q,c}.
$$

\begin{rmq}
Il est clair que chaque fonction $\dfrac{e_{q,cd}}{e_{q,c}e_{q,c}}$
est une $q$-constante, donc un \'el\'ement de $\M(\Eq)$. En
fait, on d\'eduit sans difficult\'e de la th\'eorie des
fonctions elliptiques que les $\dfrac{e_{q,cd}}{e_{q,c}e_{q,c}}$
engendrent l'extension $\M(\Eq)$ de $\C$, donc que l'on
h\'erite d'un ``gros'' corps des constantes uniquement pour
r\'esoudre les \'equations \'el\'ementaires $\sq f = c f$.
Ce n'est pas d\^u \`a un choix maladroit de fonctions de base:
quelque soit le choix des dolutions $f_{c}$ de $\sq f = c f$,
on peut d\'emontrer que les $\dfrac{f_{cd}}{f_{c}f_{c}}$ engendrent
une extension transcendante de $\C$.
\end{rmq}

% 1.4

\subsection{Le polygone de Newton} 
\label{subsection:polygone}

Toutes les preuves des assertions donn\'ees ici 
figurent\footnote{Attention: il y a eu un changement 
de convention dans la d\'efinition des pentes: celles 
\'etudi\'ees ici sont les \emph{oppos\'ees} de celles 
\'etudi\'ees dans \emph{loc.cit.}.}
dans \cite{JSFIL}. \\

Soit $P = \sum a_{i} \; \sq^{i} \in \D_{q,K}$ un op\'erateur 
aux $q$-diff\'erences non nul. On d\'efinit son \emph{polygone 
de Newton} $N(P)$ comme l'enveloppe convexe de l'ensemble:
$$
\{(i,j) \in \Z^{2} \tq j \geq v_{0}(a_{i}) \} \subset \R^{2}.
$$
C'est aussi, par d\'efinition, le polygone de Newton de l'\'equation 
aux $q$-diff\'erences $P.f = \sum a_{i} \; \sq^{i}f = 0$. 
Cette d\'efinition d\'epend du choix de la valuation $v_{0}$.
Dans le cas d'\'equations \`a coefficients dans $K_{\ell}$ (obtenues par
ramification) c'est la valuation $z_{\ell}$-adique qui sera employ\'ee.
Si l'on multiplie $P$ par un inversible $a \sq^{m}$, $a \in K^{*}$,
$m \in \Z$, le polygone de Newton est translat\'e par le vecteur
$(m,v_{0}(a)) \in \Z^{2}$. Nous consid\`ererons donc $N(P)$ comme 
d\'efini \`a une telle translation pr\`es. \\

La fronti\`ere de $N(P)$ est form\'ee de deux demi-droites verticales
et de $k \geq 1$ vecteurs de coordonn\'ees 
$(r_{1},d_{1}),\ldots,(r_{k},d_{k}) \in \N^{*} \times \Z$,
et de pentes 
$\mu_{1} = \frac{d_{1}}{r_{1}},\ldots,\mu_{k} = \frac{d_{k}}{r_{k}} \in \Q$. 
On suppose celles-ci rang\'ees par ordre croissant:
$\mu_{1} < \cdots < \mu_{k}$. La \emph{derni\`ere pente} est $\mu_{k}$
(donc, la plus grande). On notera $S(P) = \{\mu_{1},\ldots,\mu_{k}\}$ 
l'ensemble des pentes de $P$. La \emph{fonction de Newton} de $P$ est 
la fonction $r_{P}: \Q \rightarrow \N$ de support $S(P)$ et telle que 
$\mu_{i} \mapsto r_{i}$ pour $i = 1,\ldots,n$. On a donc:
$$
r_{P} = \sum_{i = 1}^{k} r_{i} \delta_{\mu_{i}},
$$
o\`u $\delta_{\mu}$ d\'esigne la fonction de Kronecker (indicatrice
de $\{\mu\}$). La correspondance entre fonctions de Newton et polygones 
de Newton (d\'efinis \`a translation pr\`es) est une bijection additive. \\

Par \emph{ramification} $z = z_{\ell}^{\ell} \;,\; q = q_{\ell}^{\ell}$, 
les pentes sont multipli\'ees par $\ell$. En particulier, en prenant 
pour $\ell$ un multiple commun des $r_{i}$, on se ram\`ene au cas o\`u 
les pentes sont enti\`eres. Cette op\'eration revient \`a une extension 
de corps aux $q$-diff\'erences, soit encore \`a une extension du corps 
de base de l'alg\`ebre $\D_{q,K}$. Elle est donc compatible avec les 
op\'erations de $\D_{q,K}$. \\

Tout \'el\'ement $\alpha$ de $K^{*}$ s'\'ecrit de mani\`ere unique
$\alpha = c z^{\mu} \beta$, o\`u $c \in \C^{*}$, $\mu \in \Z$ et 
$\beta(0) = 1$. Pour r\'esoudre l'\'equation $\sq u = \alpha u$,
nous prendrons $u = e_{q,c} \Theta_{q}^{\mu} v$, o\`u
$v(z) = \underset{k \geq 1}{\prod} \beta(q^{-k}z)$;
on v\'erifie en effet facilement que $v \in K^{*}$
(aussi bien dans le cas formel que dans le cas convergent).
On notera $e_{q,\alpha}$ l'\'el\'ement $u$ ainsi obtenu; 
c'est un \'el\'ement inversible de $L$. \\

Soit $\alpha \in K^{*}$ et soit $u$ un \'el\'ement inversible d'une
extension $K'$ de $K$ tel que $\sq u = \alpha u$. La ``transformation
de jauge'' (changement de fonction inconnue) $f = u g$ donne lieu
\`a une \'equivalence logique (changement d'\'equation)
$P f = 0 \Leftrightarrow P^{[u]} g = 0$, avec:
$P^{[u]} \underset{def}{=}  u^{-1} P u$, que l'on va expliciter.
Soit, pour simplifier, $P := \sum\limits_{i=0}^{n} a_{i} \sq^{i}$.
$$
P^{[u]} = \sum_{i=0}^{n} a_{i} u^{-1} \sq^{i} u =
\sum_{i=0}^{n} a_{i} \frac{\sq^{i}(u)}{u} \sq ^{i} =
\sum_{i=0}^{n} a_{i} \left(\prod_{j=0}^{i-1} \sq^{j}(\alpha)\right) \sq ^{i}.
$$
Ainsi, $P^{[u]} = \sum\limits_{i=0}^{n} b_{i} \sq^{i}$, avec
$b_{i} := a_{i} \prod\limits_{j=0}^{i-1} \sq^{j}(\alpha)$, donc, 
$v_{0}(b_{i}) = v_{0}(a_{i}) + i v_{0}(\alpha)$ (puisque
$v_{0}(\sq^{j}(\alpha)) = v_{0}(\alpha)$), d'o\`u l'on d\'eduit
que les pentes de $P^{[u]}$ sont 
$\mu_{1} + v_{0}(\alpha),\ldots,\mu_{k} + v_{0}(\alpha)$. On peut 
ainsi ramener une pente enti\`ere \`a $0$: si $\mu_{i} \in \Z$, 
on prend $u = e_{q,z^{-\mu_{i}}} = \Theta_{q}^{-\mu_{i}}$. \\

Remarquons que $P^{[u]}$ peut \^etre d\'efini directement en fonction 
de $\alpha$ seul, sans faire r\'ef\'erence \`a une solution de
l'\'equation $\sq u = \alpha u$ dans une extension de $K$: il suffit 
de prendre la formule \'etablie ci-dessus comme d\'efinition. 
On prouve alors facilement (exercice pour le lecteur courageux) 
que $P \mapsto P^{[u]}$ est un automorphisme de $\D_{q,K}$ et que, 
si $u,v \in K'$ sont ainsi associ\'es \`a $\alpha,\beta \in K^{*}$, 
on a la formule $P^{[uv]} = \left(P^{[u]}\right)^{[v]}$. \\

Supposons maintenant que $S(P) \subset \Z$. On va d\'efinir
l'\'equation caract\'eristique et les exposants attach\'es
\`a la $i$-\`eme pente $\mu = \mu_{i}$ de $P$. Il existe des 
indices $\alpha < \beta$ et un entier $\ell \in \Z$ tels que, 
notant $P' := P^{[e_{q,z^{-\mu}}]} = \sum a'_{i} \sq^{i}$:
\begin{align*}
\forall i \;,\; & v_{0}(a'_{i}) \geq \ell, \\
v_{0}(a'_{\alpha}) = & v_{0}(a'_{\beta}) = \ell, \\
\forall i < \alpha \;,\; & v_{0}(a'_{i}) > \ell, \\
\forall i > \beta \;,\; & v_{0}(a'_{i}) > \ell.
\end{align*}
Bien entendu, $\ell = v_{0}(P')$. Avec les notations pr\'ec\'edentes
$\mu_{i} = \dfrac{d_{i}}{r_{i}}$, on a $r_{i} = \beta - \alpha$ et
$d_{i} = v_{0}(a_{\beta}) - v_{0}(a_{\alpha})$. On introduit: 
$Q := z^{-\ell} P' = \sum b_{i} \sq ^{i}$, dont les 
coefficients sont donc dans l'anneau de valuation $\O$ de $K$
et m\^eme tels que $v_{0}(b_{\beta}) = v_{0}(b_{\alpha}) = 0$. En posant:
\begin{align*} 
\overline{Q} & \underset{def}{=} \sum b_{i}(0) \; s^{i} \\
& = b_{\alpha}(0) \; s^{\alpha} + \cdots + b_{\beta}(0) \; s^{\beta}
\in \C[s,s^{-1}],
\end{align*}
on a $b_{\alpha}(0) b_{\beta}(0) \not= 0$. L'\'equation $\overline{Q} = 0$ 
(ainsi que le polyn\^ome $\overline{Q}$ lui-m\^eme) est appel\'ee 
\emph{\'equation caract\'eristique} attach\'ee \`a la pente $\mu$ 
de $P$; on peut la consid\'erer comme d\'efinie \`a un facteur 
$c s^{k}$ pr\`es, $c \in \C^{*}, k \in \Z$. On la notera
$\overline{P}^{(\mu)}$, ou simplement $\overline{P}$ dans le cas de 
la pente $\mu = 0$. Si $\mu \not\in S(P)$, l'\'equation caract\'eristique 
est une constante non nulle. En g\'en\'eral:
$$
\overline{P}^{(\mu)} = \left(z^{-v_{0}(P^{[e_{q,z^{-\mu}}]})} 
P^{[e_{q,z^{-\mu}}]}\right)_{z = 0}.
$$
Dans tous les cas, $r_{P}(\mu)$ est \'egal au degr\'e absolu 
$\deg(\overline{P}^{(\mu)})$ de l'\'equation caract\'eristique.
L'\'equation caract\'eristique est multiplicative: 
$$
\forall P_{1}, P_{2} \in \D_{q,K}
\;,\; \forall \mu \in \Q \;,\;
\overline{P_{1} P_{2}}^{(\mu)} = 
\overline{P_{1}}^{(\mu)} \overline{P_{2}}^{(\mu)}.
$$
On en d\'eduit que le polygone de Newton est additif. Pr\'ecis\'ement, 
$P_{1}$ et $P_{2}$ \'etant des op\'erateurs aux $q$-diff\'erences 
comme ci-dessus: 
\begin{eqnarray*}
r_{P_{1}P_{2}} & = & r_{P_{1}} + r_{P_{2}}, \\
N(P_{1} P_{2}) & = & N(P_{1}) + N(P_{2}).
\end{eqnarray*}
(Cette derni\`ere \'egalit\'e concerne des parties de $\R^{2}$
d\'efinies \`a une translation de $\Z^{2}$ pr\`es !) \\

Les racines non nulles de l'\'equation caract\'eristique 
attach\'ee \`a la pente $\mu$ de $P$ sont appel\'ees les 
\emph{exposants} attach\'es \`a cette pente. Tout exposant $c$
s'\'ecrit de mani\`ere unique:
$$
c = q^{\epsilon(c)} \overline{c} \quad \text{ avec} \quad
\epsilon(c) \in \Z \quad \text{ et } 
\quad 1 \leq \lmod \overline{c} \rmod < \lmod q \rmod.
$$
Nous dirons que l' exposant $c$ attach\'e \`a la pente $\mu$ 
est \emph{non r\'esonnant} si $\epsilon(c)$ est maximal pour 
sa classe de congruence, autrement dit, si aucun $c q^{\ell}$, 
$\ell \in \N^{*}$, n'est un exposant (attach\'e \`a cette pente). \\

Si $\sq u = z^{\ell} u$, l'\'equation caract\'eristique attach\'ee 
\`a la pente $\mu$ de $P$ est \'egale \`a l'\'equation caract\'eristique 
attach\'ee \`a la pente $\mu + \ell$ de $P^{[u]}$. Si 
$\sq u = c u \;,\; c \in \C^{*}$, et si l'on note $\overline{Q}(s)$ 
l'\'equation caract\'eristique attach\'ee \`a la pente $\mu$ de $P$, 
l'\'equation caract\'eristique attach\'ee \`a la pente $\mu$ de $P^{[u]}$ 
est $\overline{Q}(c s)$. On peut donc toujours ramener un exposant 
donn\'e \`a $1$ par une telle transformation de jauge. Dans ce dernier 
cas, quitte \`a utiliser encore une transformation de jauge avec 
$u = z^{\ell} \;,\; \ell \in \Z$, on peut m\^eme supposer $1$ non r\'esonnant. 

%%%%%%%%%%%%%%%%%%%%%%%%%%%%%%%%%%%%%%%%%%%%%%%%%%%%%%%%%%%%%%%%%%%%%%%%%%%

% 2

\section{Factorisation formelle et factorisation convergente}
\label{section:factorisation}

En principe, la r\'esolution de l'\'equation (\ref{equation:eqd}) 
et la factorisation de l'op\'erateur aux $q$-diff\'erences $P$ sont 
\'etroitement imbriqu\'ees. Nous exposons d'abord la factorisation.
Pour la commodit\'e du lecteur, nous reprenons bri\`evement certains 
calculs de \cite{JSFIL}.

% 2.1

\subsection{Facteur droit associ\'e \`a un exposant non r\'esonnant}
\label{subsection:facteurdroit}

Comme on l'a vu, on peut ramener, par l'interm\'ediaire
de transformations de jauge simples, toute pente $\mu$ \`a $0$ 
et tout exposant $c$ attach\'e \`a cette pente \`a $1$. On peut 
m\^eme supposer que $1$ est non r\'esonnant, autrement dit, qu'aucun 
$q^{k} \;,\; k \in \N^{*}$ n'est un exposant attach\'e \`a la pente $0$.

\begin{lemma} \label{lemma:facteur-formel-0-1}
Supposons que $0$ est une pente de $P$ et que $1$ est un exposant 
non r\'esonnant attach\'e \`a cette pente. L'\'equation (\ref{equation:eqd})
admet alors une unique solution s\'erie formelle $f$ telle que $f(0) = 1$.
\end{lemma}
\Pr
Comme pr\'ec\'edemment (d\'efinition de l'\'equation caract\'eristique
au \ref{subsection:polygone}), on peut supposer $N(P)$ cal\'e de sorte
que la pente nulle soit sur l'axe des abcisses. On suppose de plus $P$
entier: $P = \sum\limits_{i=0}^{n} a_{i} \sq^{i}$, avec $a_{0} a_{n} \neq 0$. 
Les coefficients $a_{i}$ admettent donc un d\'eveloppement en s\'erie:
$$
\forall i \in \{0,\ldots,n\} \;,\; 
a_{i} = \sum_{j \geq 0} a_{i,j} z^{j}.
$$
L'\'equation caract\'eristique attach\'ee \`a la pente $0$ est 
donc: $\overline{P} = a_{n,0} \; s^{n} + \cdots + a_{0,0}$, et il y a
m\^eme, par hypoth\`ese, deux indices $\alpha < \beta$ tels que
$\overline{P} = a_{\beta,0} \; s^{\beta} + \cdots + a_{\alpha,0} \; s^{\alpha}$.
On \'ecrit $f = \underset{m \geq 0}{\sum} f_{m} z^{m}$ la fonction inconnue. 
Alors:
$$
P.f = \sum_{l \geq 0} g_{l} z^{l},
$$
o\`u l'on a pos\'e:
$$
g_{l} = \left(\sum_{m = 0}^{l} F_{l-m}(q^{m})f_{m}\right)
\quad , \quad
F_{j}(X) = \sum_{i = 0}^{n} a_{i,j} X^{i}.
$$
Ce polyn\^ome est constitu\'e des coefficients qui contribuent
\`a la tranche d'ordonn\'ee $j$ dans l'int\'erieur du polygone de Newton.
Ainsi $F_{0} = \overline{P}$, d'o\`u $F_{0}(1) = 0$ et
$F_{0}(q^{m}) \not= 0$ pour $m \geq 1$ (puisque $1$ est
exposant non r\'esonnant). On trouve les coefficients par
r\'ecurrence en identifiant les $g_{l}$ \`a $0$. Comme $F_{0}(1) = 0$,
$f_{0}$ est arbitraire et, comme $F_{0}(q^{m}) \not= 0$ 
pour $m \geq 1$, les $f_{m}$ sont d\'etermin\'es inductivement
de mani\`ere unique. 
\hfill $\Box$ 

\begin{lemma} \label{lemma:facteur-formel-0-c}
Supposons que $0$ est une pente de $P$. Soit $c$ un exposant 
non r\'esonnant de multiplicit\'e $m \geq 1$ attach\'e \`a la pente $0$. 
On a alors une factorisation:
$$
P = Q.(\sq  - c).u_{m}^{-1} \cdots (\sq  - c).u_{1}^{-1},
$$
o\`u $u_{1},\ldots,u_{m}$ sont des s\'eries formelles telles
que $u_{i}(0) = 1$. De plus, le polygone de Newton de $Q$ 
s'obtient en diminuant de $m$ la longueur de la pente 
horizontale de celui de $P$: $r_{Q} = r_{P} - m \delta_{0}$,
et les \'equations caract\'eristiques correspondantes 
v\'erifient: $\overline{P} = (X - c)^{m} \overline{Q}$
\end{lemma}
\Pr
Traitons d'abord le cas o\`u $c = 1$. 
Soit $u_{1} = f$, la s\'erie formelle obtenue au lemme 
\ref{lemma:facteur-formel-0-1}..
L'op\'erateur aux $q$-diff\'erences $P$ admet une factorisation:
$$
P = P_{1}.(\sq  - 1).u_{1}^{-1},
$$
o\`u l'on a pos\'e:
$$
P_{1} = \sum_{j=0}^{n-1} 
        \left(- \sum_{i=0}^{j} a_{i} \; \sq^{i}(u_{1}) \right)
        \sq ^{j}.
$$
Le calcul justificatif figure dans \cite{JSFIL} (mais il faut
modifier certains indices par suite du changement de convention
sur les pentes). Le polygone de Newton de $P_{1}$ s'obtient en 
diminuant de $1$ la longueur de la pente horizontale de celui de 
$P$: $r_{Q} = r_{P} - \delta_{0}$, et les \'equations caract\'eristiques 
correspondantes v\'erifient: $\overline{P} = (X - 1) \overline{P_{1}}$.
Il suffit alors d'it\'erer le processus pour obtenir la factorisation:
$$
P = Q.(\sq  - 1).u_{m}^{-1} \cdots (\sq  - 1).u_{1}^{-1},
$$
o\`u $u_{1},\ldots,u_{m}$ sont des s\'eries formelles telles
que $u_{i}(0) = 1$. De plus, le polygone de Newton de $Q$ 
s'obtient en diminuant de $m$ la longueur de la pente horizontale 
de celui de $P$, autrement dit, $r_{Q} = r_{P} - m \delta_{0}$. 
Les \'equations caract\'eristiques correspondantes v\'erifient: 
$\overline{P} = (X - 1)^{m} \overline{Q}$. \\

Dans le cas d'un exposant $c$ quelconque (non r\'esonnant),
soit $u = e_{q,c}$ (cf. \S \ref{subsection:polygone}). 
Alors $P^{[u]}$ v\'erifie les hypoth\`eses 
du premier cas. On \'ecrit donc:
$$
P^{[u]} = 
R.(\sq  - 1).u_{m}^{-1} \cdots (\sq  - 1).u_{1}^{-1}.
$$
On applique \`a cette \'egalit\'e la transformation de
jauge de symbole $u^{-1}$, qui commute au produit, qui
n'affecte pas les fonctions $u_{i}$ et qui transforme
$\sq  - 1$ en $c^{-1} \sq  - 1$. 
On obtient alors la factorisation voulue en prenant 
$Q = c^{-m} R^{[u^{-1}]}$. Le reste suit.
\hfill $\Box$ \\

\begin{prop} 
\label{prop:facteur-formel-mu-c}
Soit $\mu$ une pente enti\`ere de $P$. Soit $c$ 
un exposant non r\'esonnant de multiplicit\'e $m$ attach\'e 
\`a la pente $\mu$. On a alors une factorisation:
$$
P = Q.(z^{\mu} \sq  - c).u_{m}^{-1} \cdots
      (z^{\mu} \sq  - c).u_{1}^{-1},
$$
o\`u $u_{1},\ldots,u_{m}$ sont des s\'eries formelles telles
que $u_{i}(0) = 1$. De plus, le polygone de Newton de $Q$ 
s'obtient en diminuant de $m$ la longueur
de la pente de valeur $\mu$ de celui de $P$:
$r_{Q} = r_{P} - m \delta_{\mu}$,
et les \'equations caract\'eristiques correspondantes 
v\'erifient:
$\overline{P}^{(\mu)} = (X - c)^{m} \overline{Q}^{(\mu)}$.
\end{prop}
\Pr
Soit $u = e_{q,z^{-\mu}}$ (cf. \S \ref{subsection:polygone}). 
Alors $P^{[u]}$ v\'erifie les hypoth\`eses 
du lemme \ref{lemma:facteur-formel-0-c}. 
On \'ecrit donc:
$$
P^{[u]} = 
R.(\sq  - c).u_{m}^{-1} \cdots (\sq  - c).u_{1}^{-1}.
$$
On applique \`a cette \'egalit\'e la transformation de
jauge de symbole $u^{-1}$, et l'on invoque le \S \ref{subsection:polygone}.
\hfill $\Box$

\subsection{Factorisation formelle 
d'un op\'erateur aux $q$-diff\'erences}

Il y a divers \'enonc\'es possibles, en voici un:

\begin{prop} \label{prop:factorisation-formelle-souple}
Soit $\mu$ une pente enti\`ere de $P$. Soient 
$c_{1},\ldots,c_{p}$ les exposants attach\'es \`a 
la pente $\mu$, et $m_{1},\ldots,m_{p}$ leurs multiplicit\'es
respectives. On suppose les $c_{i}$ index\'es de telle
sorte que, si 
$\frac{c_{j}}{c_{i}} = q^{l} \;,\; l \in \N^{*}$,
alors $i < j$ (les exposants les moins r\'esonnants
sont factoris\'es \`a droite les premiers). On a alors une 
factorisation $P = Q R$, o\`u $\mu \not\in S(Q)$ et o\`u
$R = R_{1} \cdots R_{p}$, avec:
$$
\forall i \in \{1,\ldots,p\} \;,\;
R_{i} = (z^{\mu} \sq  - c_{i}).u_{i,m_{i}}^{-1} \cdots
        (z^{\mu} \sq  - c_{i}).u_{i,1}^{-1},
$$
les $u_{i,j}$ \'etant des s\'eries formelles telles que
$u_{i,j}(0) = 1$. On, de plus, 
$r_{R} = (m_{1} + \cdots + m_{p}) \delta_{\mu}$
et 
$\overline{R}^{(\mu)} = 
\underset{i=1}{\overset{p}{\prod}} (X - c_{i})^{m_{i}}$.
\end{prop}

\Pr
Il suffit d'appliquer r\'ep\'etitivement la proposition 
\ref{prop:facteur-formel-mu-c}.
\hfill $\Box$ \\

On peut donner des conditions plus souples. Pour la r\'esolution 
(formelle ou convergente), il sera au contraire commode d'\^etre 
plus rigide et de consid\'erer une classe d'exposants \`a la fois 
(modulo $q^{\Z}$). On obtient de la m\^eme mani\`ere:

\begin{prop} \label{prop:factorisation-formelle-rigide}
Soit $\mu$ une pente enti\`ere de $P$. Soient 
$c_{1},\ldots,c_{p}$ les exposants d'une m\^eme classe
modulo $q^{\Z}$ attach\'es \`a la pente $\mu$, 
et $m_{1},\ldots,m_{p}$ leurs multiplicit\'es respectives. 
On suppose les $c_{i}$ index\'es de telle sorte que
$\epsilon(c_{1}) < \cdots < \epsilon(c_{p})$ 
(les exposants les moins r\'esonnants sont factoris\'es 
les premiers). On a alors une unique factorisation
$P = Q R$, o\`u aucun exposant attach\'e \`a la pente
$\mu$ de $Q$ n'est congru aux $c_{i}$ modulo $q^{\Z}$
et o\`u $R = R_{1} \cdots R_{p}$, avec:
$$
\forall i \in \{1,\ldots,p\} \;,\;
R_{i} = (z^{\mu} \sq  - c_{i}).u_{i,m_{i}}^{-1} \cdots
        (z^{\mu} \sq  - c_{i}).u_{i,1}^{-1},
$$
les $u_{i,j}$ \'etant des s\'eries formelles telles que
$u_{i,j}(0) = 1$. On, de plus, 
$r_{R} = (m_{1} + \cdots + m_{p}) \delta_{\mu}$
et 
$\overline{R}^{(\mu)} = 
\underset{i=1}{\overset{p}{\prod}} (X - c_{i})^{m_{i}}$.
\end{prop}
\hfill $\Box$

\subsection{Factorisation convergente 
d'un op\'erateur aux $q$-diff\'erences}

Les r\'esultats qui pr\'ec\`edent, ainsi que leur application
\`a la r\'esolution formelle sont dus \`a Adams. Mais l'\'enonc\'e 
le plus caract\'eristique de la th\'eorie, le \emph{lemme d'Adams}, 
est l'existence de \emph{solutions convergentes associ\'ees \`a la 
derni\`ere pente} (\cite{Adams1}, \cite{Adams2}). Il se prouve le 
plus ais\'ement via la factorisation analytique, bien que celle-ci 
ait \'et\'e obtenue ult\'erieurement (\emph{voir} \cite{Birkhoff3}). \\

Nous prenons ici pour corps de base $K = \M(\C)$ ou $K = \C(\{z\})$. 
Nous dirons alors qu'une s\'erie (resp. une factorisation) est 
\emph{convergente} si elle d\'efinit un \'el\'ement de $K$ (resp. 
si tous les facteurs sont \`a coefficients dans $K$).

\begin{lemma} \label{lemma:facteur-convergent-0-1}
On reprend d'abord les hypoth\`eses du lemme 
\ref{lemma:facteur-formel-0-1}, en supposant de plus que $0$ est la 
derni\`ere pente, i.e. la plus grande. La s\'erie formelle $f$ obtenue 
comme solution est alors convergente.
\end{lemma}
\Pr
Nous commen\c{c}ons par le cas o\`u $K = \C(\{z\})$.
Avec les conventions du lemme \ref{lemma:facteur-formel-0-1}, 
la derni\`ere pente vaut $0$, les pr\'ec\'edentes sont $< 0$ et:
$$
\begin{cases}
\deg \overline{P} = \deg F_{0} = n \\
\forall j \in \{1,\ldots,n\} \;,\; \deg F_{j} \leq n
\end{cases}
$$
Ce qui suit est alors une application de la m\'ethode des
s\'eries majorantes. Des conditions sur les degr\'es et
du fait qu'aucun $F_{0}(q^{l}) \;,\, l \geq 1$ ne s'annule, 
on tire:
$$
\exists A > 0 \;:\; \forall l \in \N^{*} \;,\; 
        \left|F_{0}(q^{l})\right| \geq A \lmod q \rmod^{ln}.
$$
De la convergence des coefficients $a_{0},\ldots,a_{n}$, on tire:
$$
\exists B,C > 0 \;:\; 
    \forall i \in \{0,\ldots,n\}, \forall j \in \N \;,\;
       |a_{i,j}| \leq B C^{j},
$$
d'o\`u l'on d\'eduit (avec la condition sur les degr\'es):
$$
\forall j \in \N^{*}, \forall l \in \N \;,\;
    \left|F_{j}(q^{l})\right| \leq (n+1) B C^{j} \lmod q \rmod^{ln}.
$$
Il vient, pour $l \geq 1$:
$$
|f_{l}| \leq 
\frac{(n+1) B 
\left(
C \lmod q \rmod^{(l-1)n} |f_{l-1}| + \cdots + C^{l} \lmod q \rmod^{(l-l)n} |f_{0}|
\right)}
{A \lmod q \rmod^{ln}}.
$$
On pose $g_{l} = \frac{\lmod q \rmod^{ln} |f_{l}|}{C^{l}}$ 
et $D = \frac{(n+1)B}{A}$, et l'on a:
$$
\begin{cases}
g_{0} = 1 \\
\forall l \geq 1 \;,\; g_{l} \leq D(g_{0} + \cdots + g_{l-1})
\end{cases}
$$
On montre alors par r\'ecurrence que $g_{l} \leq (D+1)^{l}$,
d'o\`u:
$$
\forall l \geq 0 \;,\; 
|f_{l}| \leq \left(\frac{C(D+1)}{\lmod q \rmod^{n}}\right)^{l}.
$$
On a donc bien $f \in \C(\{z\})$. \\

Supposons maintenant que $K = \M(\C)$. Nous appliquons le principe, 
\'egalement caract\'eristique des \'equations aux $q$-diff\'erences, 
selon lequel ``l'\'equation fonctionnelle propage la m\'eromorphie''.
Nous r\'e\'ecrivons l'\'equation (\ref{equation:eqd}) sous la forme:
$$
f(z) = 
- \sum _{i=1}^{n} 
\frac{a_{n-i}(q^{-n}z)}{a_{n}(q^{-n}z)} f(q^{-i}z).
$$
Si $f$ est d\'efinie et m\'eromorphe dans un disque de
centre $0$ et de rayon $r > 0$, et qu'elle y v\'erifie
cette \'equation, la m\^eme formule permet de la prolonger
en une fonction m\'eromorphe sur le disque de centre $0$ 
et de rayon $\lmod q \rmod r$, qui y v\'erifie la m\^eme \'equation.
On obtient ainsi (puisque $\lmod q \rmod > 1$) un prolongement
\`a $\C$ tout entier. 
\hfill $\Box$

\begin{lemma} \label{lemma:facteur-convergent-0-c}
On reprend les hypoth\`eses du lemme \ref{lemma:facteur-formel-0-c}
en supposant de plus que $\mu$ est la derni\`ere pente. 
Les factorisations obtenues sont convergentes.
\end{lemma}

\Pr
En effet, seul le calcul de $P_{1}$ dans la premi\`ere
\'etape est non formel, et il est clair qu'il fournit
un polyn\^ome en $\sq $ \`a coefficients convergents. 
\hfill $\Box$ 

\begin{thm}(Birkhoff-Guenther) \label{thm:Birkhoff-Guenther}
On reprend les hypoth\`eses des propositions
\ref{prop:factorisation-formelle-souple}
et \ref{prop:factorisation-formelle-rigide}, en supposant
de plus que $\mu$ est la derni\`ere pente. Les factorisations 
obtenues sont alors convergentes.
\end{thm}

\Pr
Encore une fois, il suffit de rassembler les morceaux.
\hfill $\Box$

\begin{rmq}
Le rayon de convergence garanti par la preuve du lemme 
\ref{lemma:facteur-convergent-0-1}
est $\frac{\lmod q \rmod}{C(D+1)}$ : le num\'erateur d\'epend 
de $q$ seul, le d\'enominateur de l'\'equation seule.
Notons par ailleurs que cette preuve est le seul point
qui d\'epend de l'hypoth\`ese $\lmod q \rmod > 1$.
\end{rmq}

%%%%%%%%%%%%%%%%%%%%%%%%%%%%%%%%%%%%%%%%%%%%%%%%%%%%%%%%%%%%%%%%%%%%%%%%%%%

% 3

\section{Solutions formelles et solutions convergentes}
\label{section:resolution}

\emph{\textbf{Dans toute cette section, les pentes $\mu$ 
seront des entiers}}.
D'apr\`es ce qui pr\'ec\`ede, nous sommes conduits \`a nous 
int\'eresser \`a l'\'equation avec second membre: 
$$
z^{\mu} \sq f - c f = \phi.
$$ 
La transformation de jauge (cf le \S \ref{subsection:polygone})
$f = u g$, o\`u l'on choisit $u$ (dans le catalogue des fonctions 
de base) tel que $\sq u = c z^{-\mu} u$, nous ram\`ene \`a l'\'equation:
$$
\sq g - g = \gamma := \phi/ c u.
$$ 
Si l'on adopte l'analogie habituelle avec le cas diff\'erentiel:
$$
\frac{\sq - 1}{q - 1} \longleftrightarrow z \frac{d}{dz},
$$
on est conduit \`a consid\'erer cette derni\`ere r\'esolution
comme une \emph{$q$-int\'egration}. Comme dans le
cas diff\'erentiel, la constante $1$ n'est pas $q$-int\'egrable 
et n\'ecessite l'introduction du $q$-logarithme.

% 3.1

\subsection{$q$-int\'egration} 
\label{subsection:q-integration}

Soit $\pi_{0}$ le projecteur du $\C$ espace vectoriel $\C((z))$ qui 
associe \`a toute s\'erie de Laurent formelle son terme constant. 
Les sous-espaces vectoriels $\M(\C)$ et $\C(\{z\})$ sont stables,
d'o\`u, quelque soit le corps $K$, une d\'ecomposition:
$$
K = \C \oplus K^{\bullet},
\text{ o\`u }
K^{\bullet} = \left(\Ker \pi_{0}\right) \cap K.
$$
L'endomorphisme $\C$-lin\'eaire $\sq - 1$ de $K$ est nul sur la 
premi\`ere composante et laisse stable la seconde. 

\begin{lemma}
L'endomorphisme $\sq - 1$ induit un automorphisme de $K^{\bullet}$.
\end{lemma}

\Pr
En effet, on peut poser (dans $\C((z))^{\bullet}$):
$$
I_{q} \left(\sum_{i \not= 0} a_{i} z^{i}\right) =
\sum_{i \not= 0} \frac{a_{i}}{q^{i}-1} z^{i},
$$
d\'efinissant un inverse. Il est clair que celui-ci pr\'eserve,
le cas \'ech\'eant, la m\'eromorphie pr\`es de $0$ ou sur $\C$.
\hfill $\Box$ \\

On introduit donc maintenant un \'el\'ement $l_{q}$ de $L$
tel que $\sq l_{q} = l_{q} + 1$ (voir dans l'introduction les 
conventions g\'en\'erales). Noter d'ailleurs que, d'apr\`es
le calcul ci-dessus, on ne peut trouver un tel \'el\'ement
dans $K$. On note de plus, pour tout entier naturel $k$:
$$
l_{q}^{(k)} = \begin{pmatrix} l_{q} \\ k \end{pmatrix} =
\frac{1}{k!} \prod_{i=0}^{k-1} (l_{q} - i),
$$
et $l_{q}^{(k)} = 0$ pour $k < 0$, de sorte que 
(calcul facile):
$$
\forall k \in \Z \;,\;
\sq l_{q}^{(k)} = l_{q}^{(k)} + l_{q}^{(k-1)}.
$$

\begin{lemma}
Les $l_{q}^{(k)}$, $k \geq 0$, sont lin\'eairement 
ind\'ependants sur $K$; autrement dit, $l_{q}$ est
transcendant et :
$$
K[l_{q}] = \bigoplus_{k \geq 0} K l_{q}^{(k)}.
$$
\end{lemma}

\Pr
Soit en effet une relation:
$$
l_{q}^{(k+1)} = 
a_{0} \; l_{q}^{(0)} + \cdots + a_{k} \; l_{q}^{(k)},
\quad \text{les } a_{i} \in K,
$$
avec $k \geq 0$ le plus petit possible; il est donc en fait
$\geq 1$ puisque $l_{q} \not\in K$. En appliquant $\sq - 1$, 
\`a cette relation, on trouve:
$$
l_{q}^{(k)} \equiv (\sq a_{k} - a_{k}) l_{q}^{(k)}
\pmod{K l_{q}^{(0)} + \cdots + K l_{q}^{(k-1)}}.
$$
Par minimalit\'e, on en d\'eduit que
$\sq a_{k} - a_{k} = 1$, ce qui est impossible.
\hfill $\Box$ 

\begin{prop}
On a, pour tout entier naturel non nul $k$, 
une suite exacte:
$$
0 \rightarrow \C \rightarrow K_{k}[l_{q}]
\overset{\sq - 1}{\longrightarrow} K_{k-1}[l_{q}] 
\rightarrow 0.
$$
\end{prop}

\Pr
Ici, $K_{k}[X]$ d\'esigne l'ensemble des polyn\^omes de degr\'e 
$\leq k$. Ecrivons $f = f_{0} l_{q}^{(0)} + \cdots + f_{k} l_{q}^{(k)}$ 
et $g = g_{0} l_{q}^{(0)} + \cdots + g_{k-1} l_{q}^{(k-1)}$ des \'el\'ements 
respectifs de $K_{k}[l_{q}]$ et de $K_{k-1}[l_{q}]$. Par identification, 
l'\'equation $(\sq - 1) f = g$ \'equivaut \`a:
$$
\forall i \geq 0 \;,\; g_{i} = \sq f_{i} - f_{i} + \sq f_{i+1}.
$$
(On convient que $f_{i} = 0$ pour $i > k$ et que $g_{i} = 0$ 
pour $i > k-1$.) La r\'esoudre revient \`a r\'esoudre le syst\`eme:
$$
\begin{cases}
\sq f_{0} - f_{0} + \sq f_{1} = g_{0} \\
\vdots \\
\sq f_{i} - f_{i} + \sq f_{i+1} = g_{i} \\
\vdots \\
\sq f_{k-1} - f_{k-1} + \sq f_{k} = g_{k-1} \\
\sq f_{k} - f_{k} = 0
\end{cases}
$$
On voit, en commen\c{c}ant par le bas, que $f_{k} \in \C$ et m\^eme 
(avant-derni\`ere \'equation) que c'est n\'ecessairement $\pi_{0}(g_{k-1})$. 
On a alors la r\'esolution it\'erative:
\begin{eqnarray*}
f_{k} & = & \pi_{0}(g_{k-1}) \\
      & \vdots &             \\
f_{i} & = & \pi_{0}(g_{i-1}) + I_{q}(g_{i} - \sq f_{i+1}) \\
      & \vdots &             \\
f_{0} & = & \text{ une constante arbitraire } + 
            I_{q}(g_{0} - \sq f_{1})
\end{eqnarray*}
\hfill $\Box$ 

% 3.2

\subsection{Equations d'ordre $1$ avec second membre}
\label{subsection:second-membre}

On se restreint dor\'enavant \`a la sous-alg\`ebre $S$ de $L$
engendr\'ee par les fonctions \'el\'ementaires:
$$
S = K [(e_{q,cz^{\mu}})_{(c,\mu) \in \C^{*} \times \Z},l_{q}].
$$
Notons provisoirement $C(S)$ l'ensemble de tous les 
$q$-caract\`eres:
$$
C(S) = \{u \in S - \{0\} \tq
         \exists c \in \C^{*} \;:\; \sq u = c u\}.
$$
On a une d\'ecomposition
\footnote{Cette d\'ecomposition poss\`ede d'int\'eressantes 
propri\'et\'es alg\'ebriques, partiellement abord\'ees dans
\cite{SVdP} (cas formel) et \cite{JSAIF} (cas convergent).}
:
$$
S = \sum_{\mu \in \Z} S_{\mu}, \quad \text{o\`u} \quad
S_{\mu} = \sum_{u \in C(S)} u \Theta_{q}^{-\mu} K[l_{q}].
$$

La formule, imm\'ediatement v\'erifi\'ee:
$$
\sq u = c u \Rightarrow
(d z^{\nu} \sq  -1)(u \Theta_{q}^{-\mu} F) =
u \Theta_{q}^{-\mu} (cd z^{\nu - \mu} \sq  -1) F
$$
implique que l'endomorphisme $\Phi_{d,\nu} = d z^{\nu} \sq  -1$ 
du $\C$-espace vectoriel $S$ laisse stable chaque sous-espace 
$u \Theta_{q}^{-\mu} K[l_{q}]$ (et donc, en particulier, $K[l_{q}]$. 
De plus,  l'isomorphisme $F \mapsto u \Theta_{q}^{-\mu} F$ de $K[l_{q}]$ 
dans $u \Theta_{q}^{-\mu} K[l_{q}]$ conjugue l'action de $\Phi_{cd,\nu - \mu}$ 
sur le premier avec l'action de $\Phi_{d,\nu}$ sur le deuxi\`eme.
Notre but, dans ce paragraphe, est de pr\'eciser l'image et le noyau 
de ces endomorphismes. 

\begin{lemma} \label{lemma:second-membre}
Soit $(c,\mu) \in \C^{*} \times \Z$.
Il est clair que $K$ est stable par $\Phi_{c,\mu}$. \\
(i) Si $(\overline{c},\mu) \not= (1,0)$, la restriction
de $\Phi_{c,\mu}$ \`a $K$ est injective. \\
(ii) Elle est de plus surjective dans chacun des cas
suivants:
\begin{enumerate}
\item{$\mu = 0$ et $\overline{c} \not= 1$.}
\item{$\mu < 0$.}
\item{$\mu > 0$ et $K = \C((z))$.}
\end{enumerate}
\end{lemma}

\Pr
Soit $\mu = 0$. Nous \'ecrirons
$f = \underset{k > > -\infty}{\sum} f_{k} z^{k}$ et
$g = \underset{k > > -\infty}{\sum} g_{k} z^{k}$
pour rappeler que les s\'eries $f$ et $g$ n'ont qu'un nombre
fini de termes non nuls d'indice n\'egatif. On a l'\'equivalence:
$$
(c \sq  -1)f = g \Leftrightarrow
\forall k \in \Z \;,\; (c q^{k} -1) f_{k} = g_{k},
$$
qui suffit \`a montrer (i) et (ii) dans ce cas (c'est l'hypoth\`ese 
$\overline{c} \not= 1$ qui garantit que $c q^{k} - 1$ ne s'annule pas). \\

Si $\mu \not= 0$, posons $\mu = m \epsilon$, avec $m = |\mu|$ et 
$\epsilon = \pm 1$. La d\'ecomposition:
$$
\C((z)) = \bigoplus_{0 \leq i < m} z^{i} \C((z^{m}))
$$
induit des d\'ecompositions similaires de $\C(\{z\})$ et de $\M(\C)$,
et l'on \'ecrira, dans tous les cas:
$$
K = \bigoplus_{0 \leq i < m} K_{i,m}.
$$
La formule (facile \`a v\'erifier):
$$
(c z^{\mu} \sq  -1) z^{i} F(z^{m}) =
z^{i} (c q^{i} z^{m \epsilon} F(q^{m} z^{m}) - F(z^{m}))
$$
montre que chaque composante est stable. Ecrivant alors $Z = z^{m}$, 
$Q = q^{m}$, $C = c q^{i}$, $f(z) = z^{i} F(z^{m})$ et $g(z) = z^{i} G(z^{m})$, 
on obtient, en se restreignant \`a $K_{i,m}$, l'\'equivalence:
$$
(c z^{\mu} \sq  -1)f = g \Leftrightarrow
C Z^{\epsilon} F(QZ) - F(Z) = G(Z).
$$
Autrement dit, on s'est ramen\'e au cas o\`u $\mu = \pm 1$,
ce que l'on suppose maintenant. On reprend les notations
$f = \underset{k > > -\infty}{\sum} f_{k} z^{k}$ et
$g = \underset{k > > -\infty}{\sum} g_{k} z^{k}$. \\

Si $\mu = -1$, on obtient les \'equivalences:
\begin{eqnarray*}
(c z^{-1}\sq  -1)f = g & \Leftrightarrow &
\forall k \in \Z \;,\; 
c q^{k+1} f_{k+1} - f_{k} = g_{k} \\
                         & \Leftrightarrow &
\forall k \in \Z \;,\; 
c^{k+1} q^{k(k+1)/2} f_{k+1} - 
c^{k} q^{k(k-1)/2} f_{k} = 
c^{k} q^{k(k-1)/2} g_{k} \\
                         & \Leftrightarrow &
\forall k \in \Z \;,\; 
c^{k} q^{k(k-1)/2} f_{k} = 
\sum_{i < k} c^{i} q^{i(i-1)/2} g_{i}
\end{eqnarray*}
Ceci montre que $\Phi_{c,-1}$ est bijectif dans le cas formel. 
Dans le cas convergent, la relation:
$$
|c^{k} f_{k}| \leq \sum_{i < k} |c^{i} g_{i}|
$$
entraine que la s\'erie $f(cz)$ est domin\'ee par la s\'erie 
$\frac{g(cz)}{1-z}$, ce qui permet encore de conclure. \\

Si $\mu = +1$, on obtient les \'equivalences:
\begin{eqnarray*}
(c z \sq  -1)f = g & \Leftrightarrow &
\forall k \in \Z \;,\; 
c q^{k-1} f_{k-1} - f_{k} = g_{k} \\
                         & \Leftrightarrow &
\forall k \in \Z \;,\; 
\frac{f_{k-1}}{(c/q)^{k-1} q^{k(k-1)/2}} - 
\frac{f_{k}}{(c/q)^{k} q^{k(k+1)/2}} = 
\frac{g_{k}}{(c/q)^{k} q^{k(k+1)/2}} \\
                         & \Leftrightarrow &
\forall k \in \Z \;,\; 
\frac{f_{k}}{(c/q)^{k} q^{k(k+1)/2}} = 
- \sum_{i < k} \frac{g_{i}}{(c/q)^{i} q^{i(i+1)/2}}
\end{eqnarray*}
Ceci montre que $\Phi_{c,1}$ est bijectif dans 
le cas formel. 
\hfill $\Box$ \\

Dans le cas convergent avec $\mu > 0$, on ne peut pas en g\'en\'eral 
conclure, les coefficients $f_{k}$ pouvant \^etre tr\`es rapidement 
croissants. Par exemple, si $c = 1$ et $g = -1$, on trouve, pour $k > 0$,
$f_{k} = q^{k(k-1)/2}$. C'est un $q$-analogue de la s\'erie d'Euler, la
s\'erie de Tshakaloff. Pour un $g$ plus g\'en\'eral, il y a une condition
explicite pour l'existence d'une solution convergente. \\

On va maintenant \'etudier l'action de $\Phi_{c,\mu}$ sur $K[l_{q}]$. 
Le cas o\`u $(c,\mu) = (1,0)$ a fait l'objet du \S 
\ref{subsection:q-integration}. Le cas o\`u 
$(\overline{c},\mu) = (1,0)$ s'y ram\`ene car l'automorphisme 
$F \mapsto z^{l} F$ de $K[l_{q}]$ conjugue $\Phi_{c,\mu}$ avec
$\Phi_{q^{l} c,\mu}$. 

\begin{lemma}
On suppose $(\overline{c},\mu) \not= (1,0)$. Les conclusions sont 
les m\^emes que pr\'ec\'edemment : la restriction de $\Phi_{c,\mu}$ 
\`a $K[l_{q}]$ est injective; elle est de plus surjective, sauf dans 
le cas convergent si $\mu > 0$.
\end{lemma}

\Pr
Ecrivant $f = \underset{i \geq 0}{\sum} f^{(i)} l_{q}^{(i)}$
et $g = \underset{i \geq 0}{\sum} g^{(i)} l_{q}^{(i)}$
(qui sont des sommes finies), on obtient l'\'equivalence:
$$
(c z^{\mu} \sq  -1)f = g \Leftrightarrow
\forall i \geq 0 \;,\; 
(c z^{\mu} \sq  - 1) f^{(i)} =
g^{(i)} - c z^{\mu} \sq f^{(i+1)}.
$$
Ce syst\`eme se r\'esoud it\'erativement, en commen\c{c}ant
par la fin, \`a l'aide du lemme \ref{lemma:second-membre}.
\hfill $\Box$ 

\begin{cor}
On consid\`ere la restriction de $\Phi_{d,\nu}$ \`a 
$u \Theta_{q}^{-\mu} K[l_{q}]$, o\`u $\sq u = c u$. \\
(i) Si $cd = q^{l} , l \in \Z$,
et si $\mu = \nu$, cet endomorphisme est surjectif 
de noyau $\C u \Theta_{q}^{-\mu} z^{-l}$. \\
(ii) Si $\overline{cd} \not= 1$ et $\mu = \nu$,
ou bien si $cd$ est quelconque et $\nu < \mu$,
l'endomorphisme est bijectif. \\
(iii) M\^eme conclusion dans le cas formel si $\nu > \mu$.
\end{cor}

\Pr
C'est imm\'ediat par conjugaison (voir le d\'ebut
du \S \ref{subsection:second-membre}). 
\hfill $\Box$ \\

Nous synth\'etisons maintenant les r\'esultats les plus importants: 

\begin{thm} 
\label{thm:second-membre}
L'endomorphisme $\Phi_{d,\nu}$ de $S_{\mu}$ est surjectif si 
$\nu \leq \mu$, et aussi si $\nu > \mu$ dans le cas formel.
\end{thm}
\hfill $\Box$ 

% 3.3

\subsection{R\'esolution formelle}
\label{subsection:resolution-formelle}

\begin{defn}
Soient $f_{1},\ldots,f_{m}$ des \'el\'ements de $L$.
Leur $q$-Wronskien (ou Casoratien, ou Pochhammerien) est:
$$
W_{q}(f_{1},\ldots,f_{m}) =
\det \begin{pmatrix}
f_{1} & \ldots & f_{j} & \ldots & f_{m} \\
\vdots & \vdots & \vdots & \vdots & \vdots \\
\sq^{i} f_{1} & \ldots & \sq^{i} f_{j} & 
\ldots & \sq^{i} f_{m} \\
\vdots & \vdots & \vdots & \vdots & \vdots \\
\sq^{m-1} f_{1} & \ldots & \sq^{m-1} f_{j} & 
\ldots & \sq^{m-1} f_{m}
\end{pmatrix}
$$
\end{defn}

Rappelons (cf. l'introduction) que l'on note $C_{L} = L^{\sq}$ 
le sous-corps des constantes de $L$. Dans ces conditions, on a le: 

\begin{lemma} \label{lemma:q-wronskien}
Le $q$-Wronskien $W_{q}(f_{1},\ldots,f_{m})$ est non nul
si et seulement si les $f_{i}$ sont lin\'eairement ind\'ependants
sur $C_{L}$.
\end{lemma}

\Pr
Ce lemme est d\'emontr\'e dans \cite{LDVKatz}. La partie
``seulement si'' est d'ailleurs \'evidente.
\hfill $\Box$

\begin{lemma}
Le nombre de solutions ind\'ependantes de l'\'equation 
(\ref{equation:eqd}) ne peut exc\'eder $n$, l'ordre de 
l'\'equation.
\end{lemma}
\Pr
Soient en effet $f_{1},\ldots,f_{n+1}$ des solutions de l'\'equation 
(\ref{equation:eqd}). Les lignes 
$L_{i} = (\sq^{i} f_{1},\ldots,\sq^{i}f_{n+1})$
sont alors li\'ees par la relation lin\'eaire non triviale
$a_{n} L_{n} + \cdots + a_{0} L_{0}$, et $W_{q}(f_{1},\ldots,f_{n+1}) = 0$;
on conclut alors gr\^ace au lemme \ref{lemma:q-wronskien}.
\hfill $\Box$ 

Notre but est de construire une famille aussi grande que possible
de solutions ind\'ependantes de l'\'equation (\ref{equation:eqd}).
Voici un cas optimal:

\begin{thm}
Dans le cas formel, on peut construire $n$ solutions ind\'ependantes.
\end{thm}
\Pr
Elle se fait par r\'ecurrence sur l'ordre de l'op\'erateur $P$; 
l'algorithme correspondant est r\'ecursif. On exploite naturellement 
les r\'esultats sur la factorisation de la section 
\ref{section:factorisation} et ceux sur les \'equations du premier 
ordre avec second membre du \S \ref{subsection:q-integration}. \\

Si $n = 1$, on peut \'ecrire $P = a(z^{\mu} \sq - c) u^{-1}$,
et $u e_{q,c} \Theta_{q}^{-\mu}$ est une solution non nulle. \\

Si $P$ est d'ordre $n = m + 1 \geq 2$, on \'ecrit 
$P = a(z^{\mu} \sq - c) u^{-1} Q$, o\`u $Q$ est d'ordre $m$. Par 
hypoth\`ese de r\'ecurrence, il y a $m$ solutions ind\'ependantes 
$f_{1},\ldots,f_{m}$ de $Q$. D'apr\`es le th\'eor\`eme 
\ref{thm:second-membre}, il existe $f \in L$ tel que 
$Q f = u e_{q,c} \Theta_{q}^{-\mu}$. Il est clair que 
$f,f+f_{1},\ldots,f+f_{m}$ sont solutions de $P$. \\

Le d\'eterminant \'etant une forme multilin\'eaire altern\'ee,
le $q$-Wronskien de $f,f+f_{1},\ldots,f+f_{m}$ est \'egal \`a 
celui de $f,f_{1},\ldots,f_{m}$. On manipule les lignes de ce dernier: 
on remplace $L_{m}$ par $b_{m} L_{m} + \cdots + b_{0} L_{0}$, o\`u
$Q = b_{m} \sq^{m} + \cdots + b_{0}$, ce qui multiplie le d\'eterminant 
par $b_{m}$. Mais cette op\'eration remplace aussi la derni\`ere ligne 
par $(Qf,Qf_{1},\ldots,Qf_{m}) = (Qf,0,\ldots,0)$.
Le coefficient $Qf$ vaut $u e_{q,c} \Theta_{q}^{-\mu}$,
qui est inversible, et son cofacteur est le $q$-wronskien 
de $(f_{1},\ldots,f_{m})$. On obtient ainsi la formule:
$$
W_{q}(f,f+f_{1},\ldots,f+f_{m}) = 
\frac{1}{b_{m}} u e_{q,c} \Theta_{q}^{-\mu} W_{q}(f_{1},\ldots,f_{m}).
$$
Il est donc non nul, ce qui ach\`eve la preuve.
\hfill $\Box$ 

\subsection{R\'esolution analytique}

On se place ici dans le cas convergent. Si l'on reprend 
la factorisation $P = a(z^{\mu} \sq  - c) u^{-1} Q$ 
exploit\'ee au \S \ref{subsection:resolution-formelle}, 
on constate que l'on n'a la garantie
d'une factorisation convergente que si toutes les pentes
de $Q$ sont $\leq \mu$ (cf. le th\'eor\`eme
\ref{thm:Birkhoff-Guenther}). Mais, si l'une d'elles
est $< \mu$, le th\'eor\`eme \ref{thm:second-membre} ne s'applique pas.
Ainsi, la m\'ethode du \S \ref{subsection:resolution-formelle} 
ne s'applique \`a la r\'esolution convergente que si $S(P) = \{\mu\}$,
autrement dit, si $P$ est \emph{pur}.  On ne peut donc
esp\'erer trouver $n$ solutions ind\'ependantes en g\'en\'eral. \\

\begin{thm}[Lemme d'Adams]
Soit $\mu_{k}$ la derni\`ere pente de $P$. L'\'equation 
(\ref{equation:eqd}) admet alors $r_{P}(\mu_{k})$ solutions 
convergentes ind\'ependantes.
\end{thm}

\Pr
On d\'eduit en effet du \S \ref{subsection:resolution-formelle} 
une factorisation $P = QR$ avec $R$ pur de pente $\mu_{k}$ 
et d'ordre $r_{P}(\mu_{k})$. On applique alors \`a $R$ 
la m\'ethode du \S \ref{subsection:resolution-formelle} (on est dans
la cas (i) du th\'eor\`eme \ref{thm:second-membre}).
\hfill $\Box$ \\

L'exemple de l'\'equation $(\sq - 1)(z \sq - 1) f = 0$ montre 
que l'on ne peut esp\'erer mieux en g\'en\'eral.

%%%%%%%%%%%%%%%%%%%%%%%%%%%%%%%%%%%%%%%%%%%%%%%%%%%%%%%%%%%%%%%%%%%%%%%%%%%

% 4

\section{Calculs d'indices}
\label{section:indices}

Soit $P \in \D_{q,K}$. Les questions d'existence et d'unicit\'e
des solutions de l'\'equation $P.f = 0$ peuvent se traduire en
terme de surjectivit\'e et d'injectivit\'e de l'application
$f \mapsto P.f$ de $K$ dans lui-m\^eme. Nous noterons encore
$P$ cette application, qui est un endomorphisme du $\C$-espace
vectoriel $K$. \`A d\'efaut de ces propri\'et\'es optimales
(surjectivit\'e, injectivit\'e), qui \'equivalent \`a la nullit\'e
du conoyau ou du noyau de $P$, on peut tenter de calculer la
dimension de ces deux $\C$-espaces vectoriels. Nous verrons que
ces dimensions sont finies (hors les cas triviaux o\`u $P$ est
nul ou inversible). Le calcul de ces dimensions est d'ailleurs
par lui-m\^eme important pour les questions de classification
\cite{RSZ}, o\`u il est fait par d'autres m\'ethodes, combinant
l'alg\`ebre homologique et l'analyse fonctionnelle (sous la forme
de r\'esultats d\^us \`a B\'ezivin \cite{Bezivin} et \`a Ramis
\cite{RamisGrowth}).

% 4.1

\subsection{Noyau et conoyau de $\sq  - u$}

Nous commencerons par le cas d'un op\'erateur de degr\'e $1$.
\`A un facteur inversible pr\`es dans $\D_{q,K}$, on peut supposer
que $P = \sq  - u$, o\`u l'on \'ecrit $u = d q^{k} z^{\nu} v$, 
$v \in K, \; v(0) = 1$, avec $k \in \mathbf{Z}$ et 
$d \in \mathbf{C}^{*}$ tel que $1 \leq |d| < \lmod q \rmod$. 

\begin{lemma}
(i) Les dimensions du noyau et du conoyau de $\sq  - u$
ne d\'ependent que de la classe de $u \in K^{*}$ modulo
le sous-groupe $\{\frac{\sq (w)}{w} \tq w \in K^{*}\}$
de $K^{*}$. \\
(ii) Modulo ce sous-groupe, $u$ est \'equivalent \`a $d z^{\nu}$. 
\end{lemma}
\Pr
On conjugue l'endomorphisme $\mathbf{C}$-lin\'eaire 
$\sq  - u$ de $K$ par $w \in K^{*}$, identifi\'e \`a
l'automorphisme $\times w$. On trouve:
$$
w \circ (\sq  - u) \circ w^{-1} = 
\frac{w}{\sq (w)} \circ (\sq  - u'),
\text{ avec } u' = u \frac{\sq (w)}{w}.
$$
Ceci \'etablit (i). \\
On v\'erifie que $q^{k} v = \frac{\sq (w)}{w}$, o\`u
$w = z^{k} \underset{i \geq 1}{\prod} \sq ^{-i}(v) \in K$
(et que cela marche dans le cas convergent), ce qui prouve (ii).
\hfill $\Box$ \\

On peut donc supposer que $u = d z^{\nu}$, ce que nous ferons 
d\'esormais. 

\begin{lemma}
Le noyau de $\sq  - u : K \rightarrow K$ est de dimension 
$1$ si $(\nu,d) = (0,1)$, nul sinon.
\end{lemma}
\Pr
La s\'erie $f = \sum f_{n} z^{n}$ appartient au noyau si et
seulement si $\forall n \;,\, q^{n} f_{n} = d f_{n - \nu}$.
S'agissant de s\'eries de Laurent, cela entraine $f = 0$ sauf
peut-\^etre si $\nu = 0$. Dans ce dernier cas, cela entraine
$f = 0$, sauf peut-\^etre si $d \in q^{\mathbf{Z}}$. Vue la
condition sur $|d|$, ceci n'est possible que si $d = 1$, donc
$u = 1$. Dans ce cas, le noyau est $\mathbf{C}$.
\hfill $\Box$

\begin{lemma}
L'application $\sq  - u : K \rightarrow K$ est surjective
si $\nu >0$ et si $\nu = 0, \; d \not= 1$.
\end{lemma}
\Pr
Soit $g \in K$. On veut r\'esoudre (en $f$)
$\sq  f - d z^{\nu} f = g$ dans $K$. Cette \'equation 
\'equivaut \`a 
$\forall n \;,\, q^{n} f_{n} - d f_{n - \nu} = g_{n}$,
soit encore 
$\forall n \;,\, f_{n} = q^{-n} (d f_{n - \nu} + g_{n})$.
Si $\nu \geq 1$, on peut calculer les coefficients $f_{n}$
par r\'ecurrence \`a partir des $\nu$ premiers d'entre eux.
De plus, dans le cas convergent, la forme des d\'enominateurs 
montre que la s\'erie $f$ obtenue converge si $g$ converge,
et l'\'equation fonctionnelle garantit la m\'eromorphie s'il
y a lieu. Si $\nu = 0$ et $d \not= 1$ (donc en fait
$d \notin q^{\mathbf{Z}}$), on obtient directement
$f_{n} = \frac{g_{n}}{q^{n} - d}$. La convergence (resp. la
m\'eromorphie) est encore imm\'ediate le cas \'ech\'eant. \\
\emph{Variante lorsque $\nu \geq 1$:} On doit r\'esoudre 
l'\'equation au point fixe
$F(f) = f$, o\`u $F(f) := \sq ^{-1}(g + d z^{\nu} f)$.
Il est facile de voir que cet op\'erateur est contractant \`a
la fois pour la topologie $z$-adique et pour la topologie
transcendante.
\hfill $\Box$

\begin{lemma}
Si $(\nu,d) = (0,1)$, le conoyau de 
$\sq  - u : K \rightarrow K$ est de dimension $1$.
\end{lemma}
\Pr
Il est facile de voir que $\sq  - u$ annule $\mathbf{C}$
et induit un automorphisme du sous $\mathbf{C}$-espace vectoriel
$K^{\bullet}$ de $K$ form\'e des s\'eries sans terme constant.
\hfill $\Box$

\begin{lemma}
On suppose $\nu < 0$. Alors $\sq  - u : K \rightarrow K$ 
est surjectif dans le cas formel. 
\end{lemma}
\Pr
On pose $F'(f) = d^{-1} z^{- \nu} (\sq (f) - g)$. Comme
$\nu > 0$, cet op\'erateur est $z$-adiquement contractant, d'o\`u
l'existence et l'unicit\'e d'un point fixe dans le cas formel.
\hfill $\Box$ 

\begin{rmq}
C'est le premier endroit o\`u le cas formel et le cas convergent
diff\`erent. L'op\'erateur $F'$ est (fortement) \emph{dilatant}
pour la topologie transcendante, il produit des Stokes. C'est donc
ici qu'un argument analytique va \^etre n\'ecessaire dans le cas
convergent.
\end{rmq}

\begin{lemma}
On suppose $\nu = -1$. Alors, dans le cas convergent, 
le conoyau de $\sq  - u : K \rightarrow K$ est de dimension
$1$.
\end{lemma}
\Pr
Remarquons d'abord que l'automorphisme $\mathbf{C}$-lin\'eaire
$f(z) \mapsto f(dz)$ de $K$ conjugue $\sq  - d z^{-1}$
\` a $\sq  - z^{-1}$. On peut donc supposer $d = 1$ et
$u = z^{-1}$, ce que nous ferons. Par ailleurs, l'endomorphisme
$\sq  - z^{-1}$ a le m\^eme conoyau que $1 - z \sq $
et nous consid\`ererons plut\^ot ce dernier. On prendra
$K = \mathbf{C}(\{z\})$, l'autre cas s'en d\'eduisant comme 
d'habitude. \\

On consid\`ere la \emph{transformation de $q$-Borel-Ramis}
d\'efinie par:
\begin{equation}
\label{BorelRamis}
\mathcal{B}_{q,1}: 
\sum f_{n} z^{n} \mapsto \sum \frac{f_{n}}{q^{n(n-1)/2}} z^{n},
\end{equation}
qui envoie $\mathbf{C}(\{z\})$ sur l'espace 
$\mathbf{C}(\{z\})_{1}$ des s\'eries $q$-Gevrey de niveau $1$
(\emph{voir} \cite{RamisGrowth}), c'est \`a dire les
$\sum \phi_{n} z^{n}$ telles qu'il existe $A > 0$ tel que
$\phi_{n} = O(A^{n} q^{-n(n-1)/2})$ lorsque $n \to \infty$.
L'inter\^et est que la transformation $\mathcal{B}_{q,1}$ 
conjugue $1 - z \sq $ \`a la multiplication par $1 - z$,
d'o\`u le diagramme commutatif:
\begin{equation*}
\begin{CD}
\mathbf{C}(\{z\})  @>1 - z \sq >> \mathbf{C}(\{z\})  \\
@V{\mathcal{B}_{q,1}}VV                @VV{\mathcal{B}_{q,1}}V \\
\mathbf{C}(\{z\})_{1}  @>\times (1 - z)>> \mathbf{C}(\{z\})_{1}
\end{CD}
\end{equation*}
dans lequel les fl\`eches verticales sont des isomorphismes. \\
D'autre part, l'image de la fl\`eche du bas est exactement
form\'ee des s\'eries telles que $\phi(1) = 0$. En effet, l'une
des inclusions est \'evidente. Pour prouver l'autre, on prend
$\phi \in \mathbf{C}(\{z\})_{1}$ tel que $\phi(1) = 0$.
On pose $\gamma_{n} = \underset{k \leq n}{\sum} \phi_{k}$, 
de sorte que 
$\gamma(z) = \sum \gamma_{n} z^{n} \in \mathbf{C}(\{z\})$
est \'egal \`a $\frac{1}{1-z} \phi(z)$ et il reste \`a v\'erifier
la condition $q$-Gevrey sur $\gamma$. Mais la condition 
$\phi(1) = 0$ entraine 
$\gamma_{n} = -\underset{k > n}{\sum} \phi_{k}$.
Soient $C,A >0$ tels que 
$\forall n \;,\; |\phi_{n}| \leq C A^{n} |q^{-n(n-1)/2}|$.
Alors, pour $n \geq 0$:
\begin{eqnarray*}
|\gamma_{n}| & \leq & \sum_{k > n} |\phi_{k}| \\
             & \leq & C A^{n} \lmod q \rmod^{-n(n-1)/2}
               \sum_{l > 0} A^{l} \lmod q \rmod^{-l(2n+l-1)/2} \\
             & \leq & C' A^{n} \lmod q \rmod^{-n(n-1)/2}, \text{ o\`u } 
               C' = C \sum_{l > 0} A^{l} \lmod q \rmod^{-l(l-1)/2} < \infty,
\end{eqnarray*}
d'o\`u l'estimation $q$-Gevrey voulue. 
L'image de la fl\`eche du bas est donc suppl\'ementaire de
$\mathbf{C}$ et le conoyau a bien pour dimension $1$.
\hfill $\Box$ 

\begin{lemma}
On suppose $\nu = - r, r \in \mathbf{N}^{*}$. Alors, dans le cas 
convergent, le conoyau de $\sq  - u : K \rightarrow K$ 
est de dimension $r$.
\end{lemma}
\Pr
On peut \'evidemment supposer $r \geq 2$. Notons 
provisoirement $K'$ le sous-corps de $K$ form\'e des fonctions 
de $z^{r}$, d'o\`u une d\'ecomposition:
$$
K = K' \oplus z K' \oplus \cdots \oplus z^{r-1} K',
$$
dans laquelle chaque sous-espace $z^{i} K'$ est stable
par $\sq  - u$. De plus, on a un isomorphisme
$f(z) \mapsto z^{i} f(z^{r})$ de $K$ avec $z^{i} K'$,
qui conjugue $q^{i} \sigma' - d z^{-1}$ \`a $\sq  - u$,
o\`u l'on a introduit $\sigma': f(z) \mapsto f(q^{r} z)$. 
On applique alors le lemme pr\'ec\'edent.
\hfill $\Box$ \\

Nous r\'esumons ainsi nos r\'esultats:

\begin{prop} 
\label{prop:indices}
Soit $u = d z^{\nu} v$, $v \in K, \; v(0) = 1$, avec
$d \in \mathbf{C}^{*}$. Notons $\overline{d}$ la classe de $d$
modulo $q^{\mathbf{Z}}$. L'op\'erateur $\mathbf{C}$-lin\'eaire
$\sq  - u: K \rightarrow K$ est \`a indice.
Les dimensions de son noyau et de son conoyau et son indice
sont donn\'es par le tableau suivant: \\
\begin{tabular}[t]{||l|l||l|l||l||}
\hline
\multicolumn{2}{||c||}{$(\nu,\overline{d})$} & Noyau & Conoyau & 
Indice  \\
\hline
\multicolumn{2}{||c||}{$(0,1)$}   &  $1$  &  $1$ & $0$ \\
\multicolumn{2}{||c||}{$(0,\not=1)$} &  $0$  &  $0$ & $0$ \\
\multicolumn{2}{||c||}{$(> 0,-)$} &  $0$  &  $0$ & $0$ \\
\multicolumn{2}{||c||}{$(< 0,-)$} &  $0$  &  
$\begin{cases} 0 \text{ (cas formel) } \\ 
- \nu \text{ (cas convergent) } \end{cases}$ &  
$\begin{cases} 0 \text{ (cas formel) } \\ 
\nu \text{ (cas convergent) } \end{cases}$ \\
\hline
\end{tabular}
\end{prop}
\hfill $\Box$

% 4.2

\subsection{Indices et conoyaux}

Il n'y a pas en g\'en\'eral de formule simple donnant la dimension 
du noyau et du conoyau de $P$. Cependant, si toutes les pentes de
$P$ sont enti\`eres, nous savons qu'il est produit d'op\'erateurs de 
degr\'e $1$ et l'on peut d\'eduire de la proposition \ref{prop:indices} 
et de l'alg\`ebre lin\'eaire les faits suivants:
\begin{enumerate}
\item{L'application $\C$-lin\'eaire $P: f \mapsto P.f$ de $K$ dans
lui-m\^eme est \emph{\`a indice}, autrement dit, son noyau et son
conoyau sont de dimension finie.}
\item{L'\emph{indice} de $P$, c'est-\`a-dire (par d\'efinition)
l'entier $\chi(P) := \dim \Ker P - \dim \Coker P$ est la somme
des indices des facteurs de $P$.}
\end{enumerate}
En particulier, si $P$ est pur de pente $\mu \in \Z$ \emph{suppos\'ee 
non nulle}, tous ses facteurs sont de la forme $(z^{\mu} \sq -c).u$.
Dans le cas formel, ils sont tous bijectifs de $K$ dans $K$ d'apr\`es
la proposition \ref{prop:indices}, et $P$ l'est donc \'egalement.
Dans le cas convergent, chacun des $(z^{\mu} \sq -c).u$ \'etant injectif 
(toujours d'apr\`es la proposition \ref{prop:indices}), il en est de
m\^eme de $P$. Par additivit\'e de l'indice, on en d\'eduit, \emph{sans
avoir \`a supposer la pente enti\`ere}:

\begin{cor}
\label{cor:lesindices}
Soit $P$ un op\'erateur d'ordre $n$ pur de pente $\mu \neq 0$. \\
(i) Dans le cas formel, on a $\dim \Ker P = \dim \Coker P = 0$. \\
(ii) Dans le cas convergent, on a $\dim \Ker P = 0$ et
$\dim \Coker P = n \max(0,\mu)$.
\end{cor}
\Pr
Cela d\'ecoule des arguments pr\'ec\'edents si la pente est enti\`ere.
Le cas g\'en\'eral s'y ram\`ene par ramification; les d\'etails sont
laiss\'es au lecteur.
\hfill $\Box$ \\

Nous allons maintenant traduire ce r\'esultat en termes de syst\`emes,
ce qui permettra de donner une description plus concr\`ete du conoyau.
(Le but de cette article \'etant de donner des descriptions concr\`etes !)
Pour cela, nous allons d'abord d\'ecrire le lien entre \'equations et
syst\`emes aux $q$-diff\'erences. \\

Nous partons, pour simplifier, de l'op\'erateur 
$P := \sum\limits_{i=1}^{n} a_{i} \sq^{i}$, que nous supposons unitaire:
$a_{n} = 1$ (et, bien entendu, $a_{0} \neq 0$). Comme on le fait pour
les \'equations diff\'erentielles, nous allons \emph{vectorialiser}
le probl\`eme. Nous notons:
$$
X := 
\begin{pmatrix} f \\ \sq f \\ \vdots \\ \sq^{n-2} f \\ \sq^{n-1} f \end{pmatrix}
\text{~~et~~}
A := \begin{pmatrix}
0      & 1      & 0      & \ldots & 0        & 0       \\
0      & 0      & 1      & \ldots & 0        & 0       \\
\vdots & \vdots & \vdots & \ddots & \vdots   & \vdots  \\
0      & 0      & 0      & \ldots & 1        & 0       \\
0      & 0      & 0      & \ldots & 0        & 1       \\
-a_{0} & -a_{1} & -a_{2} & \ldots & -a_{n-2} & -a_{n-1}
\end{pmatrix} \in GL_{n}(K) \text{~car~} \det A = (-1)^{n} a_{0}.
$$
Il est imm\'ediat que $P.f = 0 \Leftrightarrow \sq X = A X$, 
ce qui montre comment on transforme une \'equation scalaire 
d'ordre $n$ en syst\`eme de rang $n$. Pour expliquer la r\'eciproque, 
nous devons introduire la notion de \emph{transformation de jauge}. 
Les syst\`emes de rang $n$ de matrices $A, B \in GL_{n}(K)$ sont dits 
\'equivalents si l'on peut transformer la relation $\sq X = A X$ en 
la relation $\sq Y = B Y$ par un changement de variables $Y = F X$, 
avec $F \in GL_{n}(K)$, ce qui \'equivaut \`a la relation:
$$
B = F[A] := (\sq F) A F^{-1}.
$$
(Dans $\sq F$, comme dans $\sq X$, on applique $\sq$ \`a chaque
coefficient.) Enfin, derni\`ere \'etape du raisonnement, le
\emph{lemme du vecteur cyclique}, d\^u \`a Birkhoff, dit que
tout syst\`eme de rang $n$ est \'equivalent au syst\`eme provenant
d'une \'equation d'ordre $n$ par vectorialisation; pour une preuve
\`a la main, voir \cite{JSAIF}; pour une preuve alg\'ebrique tr\`es
g\'en\'erale, voir \cite{LDVKatz}. \\

La relation entre solutions de $P.f = 0$ et solutions de $\sq X = A X$
peut \^etre pr\'ecis\'ee en une relation entre des noyaux et une
relation entre des conoyaux. Plus pr\'ecis\'ement, on a un diagramme
commutatif:
\begin{equation*}
\begin{CD}
K     @>P>>       K  \\
@V{u}VV              @V{v}VV \\
K^{n}  @>\sq - A>> K^{n}
\end{CD}
\; \text{~o\`u~}
u: f \mapsto \begin{pmatrix} 
f \\ \sq f \\ \vdots \\ \sq^{n-2} f \\ \sq^{n-1} f \end{pmatrix}
\text{~et~}
v: f \mapsto \begin{pmatrix} 0 \\ 0 \\ \vdots \\ 0 \\ f \end{pmatrix}.
\end{equation*}
La fl\`eche horizontale du bas est l'application 
$X \mapsto \sq X - A X$. Le lecteur v\'erifiera que ce diagramme
induit des isomorphismes $\Ker P \simeq \Ker(\sq - A)$ et 
$\Coker P \simeq \Coker(\sq - A)$. \\

De m\^eme, la relation entre solutions de $\sq X = A X$ et solutions 
de $\sq Y = B Y$ lorsque $B = F[A]$ peut \^etre pr\'ecis\'ee en une 
relation entre des noyaux et une relation entre des conoyaux. Plus 
pr\'ecis\'ement, on a un diagramme commutatif:
\begin{equation*}
\begin{CD}
K^{n}  @>\sq - A>> K^{n} \\
@V{F}VV           @V{\sq F}VV \\
K^{n}  @>\sq - B>> K^{n}
\end{CD},
\end{equation*}
qui fournit (encore plus facilement) des isomorphismes
$\Ker(\sq - A) \simeq \Ker(\sq - B)$ et 
$\Coker(\sq - A) \simeq \Coker(\sq - B)$. \\

Enfin, il est d\'emontr\'e dans \cite{JSAIF} que la matrice
provenant par vectorialisation d'un op\'erateur $P$ pur et
de pente $\mu$ \emph{suppos\'ee enti\`ere} est \'equivalente
\`a une matrice de la forme $A = z^{\mu} A_{0}$, o\`u 
$A_{0} \in GL_{n}(\C)$. Nous sommes donc conduits \`a \'etudier
le noyau et le conoyau de l'application 
$X \mapsto \sq X - z^{\mu} A_{0} X$ de $K^{n}$ dans lui-m\^eme. On suppose
encore $\mu \neq 0$. Il est alors imm\'ediat que cette application
est injective, dans le cas formel comme dans le cas convergent:
la relation $\sq X = z^{\mu} A_{0} X$ pour une s\'erie 
$X = \sum X_{k} z^{k}$ \'equivaut \`a la relation de r\'ecurrence
$q^{k} X_{k} = A_{0} X_{k - \mu}$; pour une s\'erie $X$ non nulle, cela 
n'est possible que si cette s\'erie admet une infinit\'e de termes
non triviaux d'exposants n\'egatifs, ce qui n'est pas autoris\'e
dans $K$. \\

\textbf{Exercice. -} 
En s'inspirant du corollaire \ref{cor:lesindices}, pr\'edire
les dimensions de $\Ker(\sq - z^{\mu} A_{0})$ et de 
$\Coker(\sq - z^{\mu} A_{0})$ dans le cas formel et dans le cas
convergent, pour $\mu > 0$ et pour $\mu < 0$. \\

Supposons $\mu > 0$. L'application $\sq - A$ est alors surjective
aussi bien dans le cas formel que dans le cas convergent. En effet:
$$
\sq X - z^{\mu} A_{0} X = Y \Longleftrightarrow
X = \sq^{-1} Y + (z/q)^{\mu} A_{0} \sq^{-1} X,
$$
que l'on peut consid\'erer comme une \'equation au point fixe, 
dont l'unique solution formelle est:
$$
X = 
\sum_{k \geq 1} q^{-\mu k(k+1)/2} z^{\mu k} A_{0}^{k} \sq^{-k-1} Y.
$$
Il n'est en effet pas difficile de voir que l'on a it\'er\'e un
op\'erateur contractant pour la distance $z$-adique, et que la
sommation ci-dessus converge formellement. Mais c'est \'egalement
un exercice \'el\'ementaire (par exemple avec des s\'eries majorantes)
que, si $Y$ converge dans un disque, alors $X$ converge dans ce disque. \\

Supposons maintenant $\mu < 0$. L'application $\sq - A$ est alors 
surjective dans le cas formel. En effet, notant $d := -\mu \in \N^{*}$:
$$
\sq X - z^{\mu} A_{0} X = Y \Longleftrightarrow
X = - z^{d} A_{0}^{-1} Y + z^{d} A_{0}^{-1} \sq X,
$$
que l'on peut encore consid\'erer comme une \'equation au point fixe, 
dont l'unique solution formelle est:
$$
X = - \sum_{k \geq 1} q^{d k(k-1)/2} z^{k d} A_{0}^{-k} \sq^{k-1} Y.
$$
En revanche, en prenant $n = 1$, $d = 1$, $Y = 1$, $A_{0} = 1$,
on voit bien que cette s\'erie n'a aucune raison de converger.
D'ailleurs, nous savons d'avance que le conoyau de $\sq - A$
est de dimension $n d$. Nous allons le v\'erifier directement.

\begin{prop}
On se place dans le cas convergent. Soient $d \in \N^{*}$ 
et $A_{0} \in GL_{n}(\C)$. L'image de l'application $\sq - z^{d} A_{0}$
de $K^{n}$ dans lui-m\^eme admet pour suppl\'ementaire $E^{n}$, o\`u 
$E := \sum\limits_{i=0}^{d-1} \C z^{i}$.
\end{prop}
\Pr
Pour $i = 0,\ldots,d-1$, notons $K_{i}$ le sous-espace de $K$ form\'e
des fonctions de la forme $z^{i} f(z^{d})$, avec $f \in K$. On a donc:
$$
K = K_{0} \oplus \cdots \oplus K_{d-1}.
$$
La relation:
$$
(\sq - z^{-d} A_{0}) \bigl(z^{i} U(z^{d})\bigr) =
z^{i} \bigl(q^{i} U(q^{d} z^{d}) - z^{-d} A_{0} U(z^{d})\bigr)
$$
montre que chaque $K_{i}$ est stable; de plus, notant $z' := z^{d}$,
$K' := \{f(z^{d}) \tq f \in K\}$, $q' := q^{d}$ et 
$\sigma': f(z') \mapsto f(q' z')$, la restriction
de $\sq - z^{-d} A_{0}$ \`a $K_{i}$ est conjugu\'ee par l'isomorphisme
$f \mapsto z^{i} f$ de $K'$ sur $K_{i}$ \`a l'application
$q^{i} \sigma' - (z')^{-1} A_{0}$ de $\C(\{z'\}^{n}$ dans lui-m\^eme. 
On est ainsi ramen\'e au cas $d = 1$. \\
Il s'agit alors de d\'emontrer que l'image de l'application 
$\sq - z^{-1} A_{0}$ de $K^{n}$ dans lui-m\^eme admet pour 
suppl\'ementaire $\C^{n}$. Il est \'equivalent, et plus commode,
de consid\'erer 'application $z \sq - A_{0}$. On \'ecrit cela sous 
la forme d'une d\'ecomposition:
$$
Y = z \sq X - A_{0} X + Z;
$$
plus pr\'ecis\'ement, il s'agit, $Y \in K^{n}$ \'etant donn\'e,
de d\'eterminer $X \in K^{n}$ et $Z \in \C^{n}$ (qui doivent en
principe \^etre uniques). On \'ecrit $X = \sum X_{k} z^{k}$,
$Y = \sum Y_{k} z^{k}$ (attention, cette notation n'a pas de rapport
avec celle des $K_{i}$ plus haut). La relation de r\'ecurrence qui 
vient est:
\begin{align*}
Y_{0} & = q^{-1} X_{-1} - A_{0} X_{0} + Z, \\
Y_{k} & = q^{k-1} X_{k-1} - A_{0} X_{k} \; (k \neq 0).
\end{align*}
En appliquant la transformation de Borel-Ramis (\ref{BorelRamis}),
page \pageref{BorelRamis}, on est ramen\'e \`a l'\'equation:
$$
\mathcal{B}_{q,1} Y = (z - A_{0}) \mathcal{B}_{q,1} X + Z.
$$
Un calcul similaire \`a celui d\'ej\`a fait montre alors que
l'unique solution s'obtient en prenant:
$$
Z = \sum_{k \in \Z} q^{-k(k-1)/2} A_{0}^{k} Y_{k}.
$$
\hfill $\Box$

%%%%%%%%%%%%%%%%%%%%%%%%%%%%%%%%%%%%%%%%%%%%%%%%%%%%%%%%%%%%%%%%%%%%%%%%%%%%%

\end{document}